\documentclass[final,3p,times]{elsarticle}
\biboptions{numbers}

\usepackage{hyperref}
\usepackage{cancel}
\usepackage{graphicx}%
\usepackage{multirow}%
\usepackage{amsmath,amssymb,amsfonts}%
\usepackage{amsthm}%
\usepackage{mathrsfs}%
\usepackage[title]{appendix}%
\usepackage{xcolor}%
\usepackage{textcomp}%
\usepackage{manyfoot}%
\usepackage{booktabs}%
\usepackage{algorithm}%
\usepackage{algorithmicx}%
\usepackage{algpseudocode}%
\usepackage{listings}%

\usepackage{amsfonts}
\usepackage{graphics}
\usepackage{bm}
\usepackage{bbm}
\usepackage[section]{placeins}
\usepackage{caption}
\usepackage{subcaption}
\usepackage{epsfig}
\usepackage{booktabs}
\usepackage{color}
\usepackage{mathtools}
\graphicspath{ {./figures/} }
\usepackage{booktabs}
\usepackage{multirow}
\usepackage{pifont}
\usepackage{inputenc}
\usepackage{tikz}
\usetikzlibrary{shapes,arrows}
\usetikzlibrary{shapes.geometric}
\usetikzlibrary{calc}

\newcommand{\ary}[1]{\boldsymbol{\mathsf{#1}}}
\newcommand{\pder}[2]{\frac{\partial #1}{\partial #2}}
\newcommand{\pdersq}[2]{\frac{\partial^2 #1}{\partial #2^2}}
\newcommand{\pderiv}[3]{\frac{\partial^2 #1}{\partial #2 \partial #3}}

\usepackage{lineno}
\usepackage{mathrsfs}
\usepackage{color}
\usepackage{xcolor}

\journal{Elsevier}

\begin{document}

\begin{frontmatter}



\title{Unveiling the Multiphysics Complexity: An Isogeometric Framework for Inducing Bifurcation and Tracing Post‑Buckling Paths in Electroelastic Thin Shells}

\author[label1]{Zhaowei Liu\corref{cor1}}
\affiliation[label1]{organization={College of Mechanics and Engineering Science},
            addressline={Hohai University},
            city={Nanjing},
            postcode={211000},
            state={Jiangsu},
            country={China}}
\ead{zhaowei.liu@hhu.edu.cn}

\author[label1]{Long Jin}

\author[label3]{Andrew McBride}
\address[label3]{Glasgow Computational Engineering Centre, James Watt School of Engineering, University of Glasgow, Glasgow, G12 8LT, United Kingdom}
\author[label2]{Weicheng Huang}
\affiliation[label2]{organization={School of Engineering},
            addressline={Newcastle University, Stephenson Building},
            city={Newcastle upon Tyne},
            postcode={NE1 7RU},
            country={United Kingdom}}
\author[label1]{Tiantang Yu}
\author[label1]{Peiliang Bian}
\author[label3,label4]{Paul Steinmann}

\address[label4]{Institute of Applied Mechanics, Friedrich-Alexander Universit\"at Erlangen-N\"urnberg, D-91052, Erlangen, Germany}

\cortext[cor1]{corresponding authors}
\begin{abstract}
Electroelastic shells are widely used in soft actuators, sensors, and energy harvesters owing to their large electrically induced deformations. 
However, the accurate simulation of their complex nonlinear multiphysics coupling, including bifurcation and post-buckling responses, remains challenging.
This work presents an isogeometric Kirchhoff--Love shell formulation for the nonlinear analysis of electroelastic thin structures undergoing finite deformations. The formulation incorporates geometrically nonlinear kinematics, Maxwell-stress-induced electromechanical coupling, material incompressibility, and initial prestretch. Catmull--Clark subdivision surfaces are employed to ensure the $C^1$ continuity required by Kirchhoff--Love shell theory. Consistent tangent operators are derived analytically, and a static condensation procedure is introduced to satisfy the plane-stress constraint. To trace bifurcation and post-buckling equilibrium paths, a staged Newton--Raphson algorithm with arc-length continuation and eigenmode perturbation is adopted. Numerical examples involving spherical membranes, prestretched circular plates, and toroidal membranes demonstrate the capability of the proposed framework to accurately capture large deformations, symmetry-breaking instabilities, and post-buckling responses under coupled electromechanical loading.
\end{abstract}




\end{frontmatter}


\section{Introduction}
\label{intro}
Dielectric elastomers (DEs) represent a class of electroactive polymers that exhibit significant deformation under electric fields, making them promising for actuators, sensors and energy harvesters. 
The theoretical foundations of electroelasticity were established by~\citet{toupin1956elastic}, who formulated a general theory for elastic dielectrics, followed by contributions from~\cite{pao1978electromagnetic,maugin2013continuum,eringen2012electrodynamics}. 
These works laid the groundwork for nonlinear continuum electromechanics, incorporating Maxwell stresses and polarisation effects. 
\citet{dorfmann2005nonlinear,dorfmann2014nonlinear} developed a comprehensive nonlinear theory for electroelasticity, introducing constitutive models based on free energy functions. 
This was extended to incompressible materials by~\citet{bustamante2009nonlinear}, who derived variational principles for coupled problems. 
The modeling of DEs often involves hyperelastic potentials, such as the neo-Hookean and Gent models, to capture large deformations~\cite{vu2007numerical,zah2015multiplicative}.

Constitutive modeling for DEs has evolved to account for material nonlinearity and electromechanical coupling. 
\citet{vu2007numerical} proposed a finite element formulation for EAPs using a free energy-based approach. 
\citet{skatulla2012multiplicative} introduced a multiplicative formulation for nonlinear electro-elasticity, while~\citet{zah2015multiplicative} incorporated micromechanically motivated network models. 
For nearly incompressible materials, \citet{dorfmann2010nonlinear} addressed the challenges of volumeric locking by using mixed formulations. 
\citet{pechstein2020large} extended these ideas to large deformations, employing augmented free energy functions.
Experimental characterisation by~\citet{hossain2012experimental} and~\citet{mehnert2021complete} provided data for model validation, highlighting the importance of accurate permittivity and hyperelastic parameters.

The development of robust finite element methods for DEs has been subject to focus due to locking issues in thin structures. 
Early work by~\citet{sze1999hybrid} adopted special stress elements for piezoelectric materials. 
\citet{klinkel2006geometrically} proposed a mixed formulation with six independent fields, and subsequently extended it to nonlinear dielectrics~\cite{klinkel2013solid}. 
The Tangential Displacement Normal Normal Stress (TDNNS) method, introduced by~\citet{pechstein2011tangential}, avoided shear and volume locking by using mixed elements with tangential displacement continuity. 
This was applied to piezoelectric solids~\cite{pechstein2018new} and later to large-deformation electro-elasticity~\cite{pechstein2020large}. 
For shells, \citet{neunteufel2019hellan} developed a Hellan-Herrmann-Johnson-type formulation, which~\citet{pechstein2025efficient} adapted to dielectric elastomer shells with independent thickness deformation.
\citet{libai2012nonlinear} and~\citet{vetyukov2014nonlinear} provided foundations for nonlinear shell theory. 
For DEs, \citet{ortigosa2016new} presented a convex multi-variable potential for large strains, while ~\citet{klinkel2013solid} developed a solid shell element with through-thickness electric field approximation. 
\citet{pechstein2020large} introduced relaxed Kirchhoff-Love kinematics with independent thickness stretch, validated against 3D benchmarks. 
Applications include buckling actuators~\cite{pechstein2020large}, peristaltic pumps~\cite{lotz2009peristaltic}, and spherical grippers~\cite{kadapa2020robust}, demonstrating the ability to capture complex instabilities.

Computational efficiency is a paramount concern in the numerical simulation of dielectric elastomers, as these materials exhibit complex behaviours like near-incompressibility and geometric nonlinearities that challenge conventional finite element methods. 
A significant contribution comes from \citet{kadapa2020robust}, who introduced a novel framework employing Bézier elements within a mixed displacement-pressure formulation. 
This approach leverages monolithic solving strategies to overcome key limitations of traditional elements, such as Q1/P0 and F-bar elements \cite{de2011computational}, which often suffer from volumetric locking and poor convergence for incompressible materials. 
The framework's effectiveness lies in its ability to maintain stability under large deformations.

Conventional finite element methods, based on Lagrange polynomials, often struggle with the Kirchhoff-Love shell formulation due to the requirement for $C^1$-continuous discretisations. 
This continuity condition ensures proper representation of bending effects without rotational degrees of freedom. 
To achieve $C^1$-continuity in Kirchhoff-Love shell discretisations, specialised approaches include exotic finite elements like Argyris spaces on triangles~\cite{Argyris1968} and quadrilaterals~\cite{Kapl2021}, TUBA plate elements for adaptive $C^1$ discretisation~\cite{ivannikov2015generalization}, as well as discontinuous Galerkin~\cite{noels2008new} and meshless methods~\cite{krysl1996analysis,ivannikov2014meshless} for weak continuity imposition. 
These methods provide robust alternatives to conventional formulations while maintaining computational efficiency.
Isogeometric analysis (IGA)~\cite{hughes2005isogeometric} has emerged as a powerful alternative, offering smooth basis functions that naturally satisfy these continuity requirements. 
Recent advances in IGA Kirchhoff–Love shell formulations have addressed membrane locking 
through computationally efficient discretisations~\cite{mathews2024computationally} and have extended the framework to trimmed multi-patch geometries using reduced-order 
methods~\cite{chasapi2024fast}. 
Among IGA approaches, subdivision surfaces~\cite{cirakortiz2000,Cirak:2001aa} provide particularly attractive features, including the ability to handle complex geometries with arbitrary topology while maintaining the desired smoothness~\cite{liu2024computational}.

However, accurately capturing symmetry-breaking and post-buckling paths induced by electromechanical coupling instabilities in the numerical simulation of dielectric elastomer shells remains a challenge~\cite{langham2018modeling,feng2021numerical,sun2022snap}. 
Recent studies have explored tunable morphing of DE balloons~\cite{su2023tunable} 
and exploited instabilities for large shape transformations in DEs~\cite{katusele2025exploiting}, 
yet a unified numerical framework capable of robustly tracing post-bifurcation paths 
in thin-shell geometries under combined electromechanical loading is still lacking.
This requires high-fidelity computational models capable of handling large deformations, and geometric and material nonlinearities. 
The IGA analysis framework based on subdivision surfaces proposed in this paper, owing to its high-order continuity and accurate geometric representation, is particularly well-suited for simulating such nonlinear phenomena involving smooth deformation modes and complex instability patterns.
This paper presents a comprehensive framework for the isogeometric analysis of electroelastic thin shells based on subdivision surfaces. 
The main contributions of this work are:
\begin{enumerate}
\item A nonlinear Kirchhoff–Love shell formulation specifically developed for dielectric elastomers, incorporating finite deformation kinematics and electromechanical coupling effects, with a consistent treatment of Maxwell stress and material incompressibility.
\item A systematic numerical framework for inducing bifurcation and tracing post-buckling equilibrium paths in electroelastic thin shells. This is achieved through a staged arc-length procedure combined with eigenmode perturbation, enabling the robust detection of symmetry-breaking instabilities and the stable traversal of unstable equilibrium branches.
\item The use of Catmull–Clark subdivision surfaces provides the $C^1$-continuity required by Kirchhoff–Love shell theory, ensuring a smooth representation of deformation fields even in the presence of severe localisation and self-contact during post-buckling.
\item Through comprehensive numerical examples, including spherical membranes, prestretched circular plates, and toroidal membranes, the proposed method is validated against analytical solutions, demonstrating its unique capability to capture bifurcation onset, mode switching, and post-buckling responses under coupled electromechanical loading.
\end{enumerate}

The remainder of this paper is organised as follows. Section~\ref{sec:dielectricity} establishes the theoretical foundation for electromechanical coupling in dielectric elastomers. 
Section~\ref{sec:elec_shells} provides a comprehensive electroelastic shell formulation, addressing geometric description, kinematics, constitutive modelling, stress decomposition, incompressibility and plane stress constraints, stress resultants, consistent tangent moduli via static condensation, and the treatment of initial prestretch. Section~\ref{sec:implementation} details the numerical implementation using subdivision surfaces. 
Section~\ref{sec:bifurcation} describes the specialised techniques and algorithms developed for analysing bifurcation and post-buckling behaviour in electroelastic shells. Section~\ref{sec:numerical_examples} provides numerical examples that validate and showcase the proposed approach, and Section~\ref{sec:conclusions} concludes with a summary of key findings and future research directions.

\section*{Notations}
\label{sec:notations}
\subsection*{Brackets}
Square brackets $[ \, ]$ are used to group algebraic expressions. Round brackets $( \,)$ are used to denote the dependencies of a function.  If brackets are used to denote an interval, then $(\,)$ stands for an open interval and $[\,]$ is a closed interval. Curly brackets $\{\,\}$ are used to define sets.
\subsection*{Symbols}
A variable typeset in a normal weight font represents a scalar. A bold weight font denotes a first- or second-order tensor. An overline indicates that the variable is defined with respect to the reference configuration. If absent, the variable is defined with respect to the deformed configuration. A scalar variable with superscript or subscript indices normally represents the components of a vector or second-order tensor. Upright font is used to denote matrices and vectors.

Indices $i,j,k,\dots$ vary from $1$ to $3$, while $\alpha, \beta, \gamma,\dots$, used to indicate surface variable components, vary from $1$ to $2$. Einstein summation convention is used throughout.

The comma symbol in a subscript represents a partial derivative, for example, $A_{,\beta}$ is the partial derivative of $A$ with respect to the $\beta^{\text{th}}$ coordinate.

{
To ensure clarity and avoid ambiguity in the nonlinear formulation and numerical implementation, the following notation is adopted throughout this work:
\begin{enumerate}
    \item \textbf{Layer Difference:} 
    $\Delta$ denotes the physical difference between the upper and lower layers of the shell 
    (e.g., $\Delta \Phi = \Phi_{\text{top}} - \Phi_{\text{bottom}}$).

    \item \textbf{Variation:} 
    $\delta$ denotes the first variation of a variable 
    (e.g., $\delta \mathbf{u}$ is the virtual displacement).

    \item \textbf{Linearisation:} 
    $\varDelta$ denotes the linearisation (total increment) over a solution step 
    (e.g., $\varDelta \mathbf{u}$ is the displacement increment for the current arc-length step).

    \item \textbf{Iterative Improvements:} 
    Symbols with a tilde ($\tilde{\cdot}$) are reserved for iterative updates within the nonlinear solver.
    Within a Newton--Raphson iteration, the iterative correction $\tilde{\delta}\mathbf{u}$ is added to the accumulated increment:
    \[
        \tilde{\varDelta}\mathbf{u}^{(k+1)}
        = \tilde{\varDelta}\mathbf{u}^{(k)} + \tilde{\delta}\mathbf{u},
    \]
    so that the total increment $\tilde{\varDelta}\mathbf{u}$ is the sum of all iterative improvements.
\end{enumerate}
}
\subsection*{Coordinates}
$x$, $y$, and $z$ denote the Cartesian coordinates of a three-dimensional Euclidean space.
${\theta}^i$ denotes coordinates in the local element space. The three covariant basis vectors for a surface point are denoted as $\mathbf a_i$, where $\mathbf a_1, \mathbf a_2$ are tangential vectors and $\mathbf a_3$ is the normal vector. 

\section{Electromechanical Coupling in Dielectric Elastomers}
\label{sec:dielectricity}
This section presents the three-dimensional electroelastic framework that serves as the foundation of the proposed shell formulation. 
The electrostatic governing equations and constitutive relations are first introduced in both the reference and current configurations. 
Subsequently, the transformation of electric field quantities under finite deformation is described, followed by the formulation of a general electroelastic strain-energy density function to characterise the coupling between mechanical deformation and electric fields.
\subsection{Configurations} 
\begin{figure}
    \centering
    \includegraphics[width=0.7\linewidth]{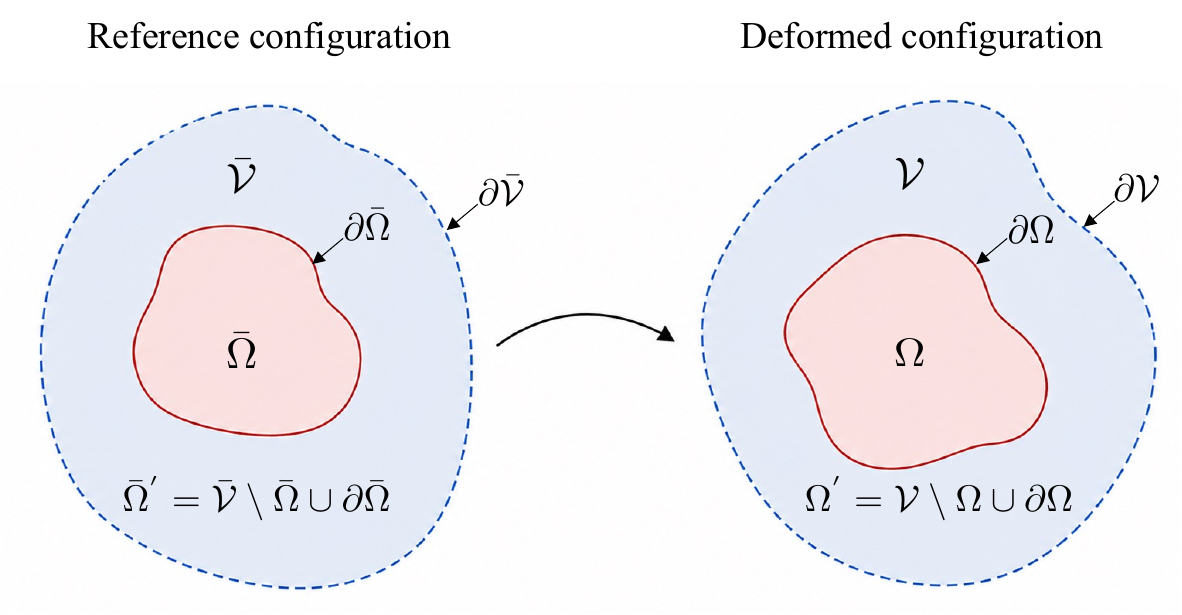}
    \caption{Reference and deformed configurations of a dielectric elastomer and the surrounding free space.}
    \label{fig:configurations_3d}
\end{figure}

To analyse the electromechanical coupling problem involving a dielectric elastomer, it is necessary to carefully consider the governing equations and the relationship between the electric field and electric displacement in reference and deformed configurations.
Consider the dielectric elastomer occupying a region of space $\bar{\Omega}$ and ${\Omega}$ in ${\mathbb{R}}^3$ in its reference and deformed configuration, respectively (shown in Fig.~\ref{fig:configurations_3d}).
One introduces the deformation map $\mathbf{r} = \chi(\mathbf{\bar{r}})$ to describe the motion from the reference to the deformed configuration, where $\mathbf{r}$ and $\mathbf{\bar{r}}$ are the position vectors in deformed and reference configurations, respectively.
Also, the three-dimensional region $\Omega$ lies inside region $\mathcal{V}$, so that
the surrounding free space is
$ {\Omega}^{\mathrm{'}} = \mathcal{V} \setminus \Omega \cup \partial \Omega.$
 For $\bar{\mathcal{V}}$ is the referential region corresponding to $ \mathcal{V}$ in ${\mathbb{R}}^3$, then
$
   \bar{\Omega}^{'}=  \bar{\mathcal{V}} \setminus \bar{\Omega} \cup \partial \bar{\Omega}.
$
\subsection{Electrostatic Governing Equations}
The electrostatic problem for a dielectric elastomer is governed by Maxwell's equations in the deformed domain $\mathcal{V}$.
The spatial electric field $\mathbf{E}$ satisfying Faraday's Law for electrostatics as
\begin{equation}
\mathrm{\nabla} \times\,\mathbf{E}=\mathbf{0}.
\end{equation}
This implies that the electric field is irrotational and it can be expressed as the gradient of a scalar potential $\Phi$ expressed as
\begin{equation}
\mathbf{E}=-\mathrm{\nabla}\,\Phi.
\label{eq:H-Phi_relation}
\end{equation}
Also, in the absence of volume charges, the spatial electric displacement $\mathbf{D}$ is governed by the Gauss's Law for Electricity as
\begin{equation}
\mathrm{\nabla} \cdot\,\mathbf{D} = 0.
\label{eq:Mag_field}
\end{equation}
\subsection{Constitutive Relationship} 
In a dielectric material, the electric displacement is related to the electric field by a constitutive equation as
\begin{equation}
\mathbf{D}=\epsilon \mathbf{E} + \mathbf{P},
\label{eq:magneticc_constitutive}
\end{equation}
where $\epsilon$ is the constant electric permittivity of free space and $\mathbf{P}$ is the spatial polarisation, which vanishes in $\Omega^{'}$. 
\subsection{Transformation Between Configurations}
To relate the electric displacement and electric field in the reference and deformed configurations, the deformation gradient tensor $\mathbf{F} := \partial\mathbf{r}/\partial\bar{\mathbf{r}}$, is adopted to perform a pull-back transformation, where $\bar{\mathbf{r}}$ and $\mathbf{r}$ are the position vectors of a material point in the reference and deformed configurations, respectively. Thus, the referential electric displacement $\bar{\mathbf{D}}$ and electric field $\bar{\mathbf{E}}$ (defined over the reference domain $\bar{\mathcal{V}}$) are given by
\begin{eqnarray} 
\bar{\mathbf{D}} &=& \mathcal{J} \mathbf{F}^{-1} \mathbf{D}, \nonumber \\
\bar{\mathbf{E}} &=& \mathbf{F}^{\mathrm{T}}\mathbf{E},
\label{eq:transformations}
\end{eqnarray}
where $\mathcal{J} = \mathrm{det}(\mathbf{F})$ represents the local volume change due to deformation. The electric field does not carry the Jacobian factor $\mathcal{J}$ under the pull-back because it is a vector quantity that transforms covariantly with the deformation gradient. In contrast, the electric displacement is a flux density and must be scaled by $\mathcal{J}$ to account for the change in cross-sectional area during deformation. These referential quantities satisfy the following Maxwell's equations in the reference configuration:
\begin{align}
\bar{\nabla}\cdot\,\bar{\mathbf{D}} &= 0, \nonumber \\
\bar{\nabla}\times\,\bar{\mathbf{E}} &= \mathbf{0},
\end{align}
where $\bar{\nabla}$ denotes the gradient operator with respect to the reference coordinates $\bar{\mathbf{r}}$. The electric field in the reference configuration can be expressed as the gradient of a scalar potential $\bar\Phi$:
\begin{equation}
\bar{\mathbf{E}} = -\bar{\nabla} \bar\Phi.
\end{equation}
The scalar potential $\bar{\Phi} = \bar{\Phi}(\bar{\mathbf{r}})$ is the referential counterpart of the spatial potential $\Phi = \Phi(\mathbf{r})$, obtained via the \emph{pull-back} operation through the deformation map $\mathbf{r} = \boldsymbol{\chi}(\bar{\mathbf{r}})$:
\begin{equation}
\bar{\Phi}(\bar{\mathbf{r}}) = \Phi\big(\boldsymbol{\chi}(\bar{\mathbf{r}})\big).
\label{eq:potential_map}
\end{equation}
This equality ensures that the potentials coincide at each material point, despite being expressed as functions of different coordinates.

Using the transformations~\eqref{eq:transformations}, the constitutive relationship~\eqref{eq:magneticc_constitutive} in the reference configuration is derived from its spatial form as
\begin{equation}
\mathcal{J}^{-1} \mathbf{C}\,\bar{\mathbf{D}} = \epsilon\,\bar{\mathbf{E}} + \bar{\mathbf{P}} \quad \text{in } \bar{\mathcal{V}},
\label{eq:magneticc_constitutive_ref}
\end{equation}
where $\bar{\mathbf{P}} = \mathbf{F}^{\mathrm{T}} \mathbf{P}$ and $\mathbf{C} = \mathbf{F}^{\mathrm{T}}\mathbf{F}$ is the right Cauchy-Green deformation tensor. Since the polarization $\mathbf{P}$ vanishes identically in the free space $\Omega^{\prime}$ surrounding the material, it follows that $\bar{\mathbf{P}} = \mathbf{F}^{\mathrm{T}}\mathbf{P} = \mathbf{0}$ in the corresponding referential region $\bar{\Omega}^{\prime}$. Consequently, the constitutive relationship simplifies to
\begin{equation}
\bar{\mathbf{D}} = \epsilon \mathcal{J} \mathbf{C}^{-1} \bar{\mathbf{E}} \quad \text{in } \bar{\Omega}^{\prime}.
\label{free_sp_constitutive}
\end{equation}
\subsection{Strain Energy Density Function}
\label{sec:strain_energy}
{
For a general incompressible electroelastic solid one may regard the strain energy density function in the reference configuration as a function of the right Cauchy--Green tensor $\mathbf{C}$ and the reference electric field vector $\bar{\mathbf{E}}$ as:

\begin{equation}
    W(\mathbf{C},\bar{\mathbf{E}}) = \widetilde{W}(I_1,I_2,I_3,I_4,I_5,I_6),
    \label{eq:energy_density}
\end{equation}
where $\widetilde{W}$ depends on the six scalar invariants of $\mathbf{C}$ and $\bar{\mathbf{E}}$. The first three invariants are purely mechanical:
\begin{align}
    I_1 &= \operatorname{tr}\mathbf{C}, &
    I_2 &= \frac{1}{2}\big[(\operatorname{tr}\mathbf{C})^2 - \operatorname{tr}(\mathbf{C}^2)\big], &
    I_3 &= \det\mathbf{C} = \mathcal{J}^2,
    \label{eq:mechanical_invariants}
\end{align}
The remaining invariants involve the electric field:
\begin{align}
    I_4 &= \bar{\mathbf{E}} \cdot \bar{\mathbf{E}}, &
    I_5 &= \bar{\mathbf{E}} \cdot [\mathbf{C}\bar{\mathbf{E}}], &
    I_6 &= \bar{\mathbf{E}} \cdot [\mathbf{C}^2\bar{\mathbf{E}}],
    \label{eq:electrical_invariants}
\end{align}
where $I_4$ captures purely electric effects, while $I_5$ and $I_6$ represent electromechanical coupling.
For incompressible materials, $\mathcal{J} \equiv 1$, and consequently $I_3 \equiv 1$. 
A common simplifying assumption~\cite{dorfmann2005nonlinear,zah2015multiplicative}, adopted in the present work, is that the mechanical and electrical contributions are \emph{additively separable}:
\begin{equation}
    \widetilde{W}= \widetilde{W}_{\mathrm{mech}} + \widetilde{W}_{\mathrm{elec}}.
    \label{eq:energy_decomposed}
\end{equation}

\paragraph{Electric energy for voltage‑controlled elastomers}
In the present work, the dielectric material is assumed to be isotropic and linearly polarizable. Under voltage control, the electric energy density~\cite{zah2015multiplicative} is most naturally written in the current configuration as
\begin{equation}
    \widetilde{W}_{\mathrm{elec}}(\mathbf{C},\bar{\mathbf{E}})
    = -\frac{1}{2}\,\epsilon\,\mathbf{E}\cdot\mathbf{E}
    = -\frac{1}{2}\,\epsilon\,
      \bigl[\bar{\mathbf{E}}\otimes\bar{\mathbf{E}}\bigr]:\mathbf{C}^{-1}.
    \label{eq:elec_energy_1}
\end{equation}
Using the Cayley–Hamilton theorem for incompressible materials, this expression can be rewritten entirely in terms of the invariants:
\begin{equation}
    \widetilde{W}_{\mathrm{elec}}
    = -\frac{1}{2}\,\epsilon\,
      \bigl[I_6 - I_1 I_5 + I_2 I_4\bigr].
    \label{eq:elec_energy_invariants}
\end{equation}
For comprehensive derivations, refer to~\ref{app:elec_energy_invariants}.}
\section{Electroelastic Shells}
\label{sec:elec_shells}
The analysis of electroelastic thin shells presents several challenges that fundamentally distinguish it from purely mechanical shell problems. First, the presence of an electric field introduces additional stress contributions, namely Maxwell stresses, which are strongly coupled with mechanical deformation and give rise to pronounced nonlinear electromechanical interactions~\cite{klinkel2006geometrically,klinkel2013solid}. Second, the thin-shell geometry necessitates a careful treatment of through-thickness kinematics and the plane stress condition, particularly under large deformations where thickness stretching becomes non-negligible. Third, the incompressibility constraint characteristic of dielectric elastomers must be incorporated consistently within the shell framework to ensure physically admissible deformation states.

To address these challenges, this section develops a comprehensive electroelastic shell formulation based on the Kirchhoff--Love hypothesis. The formulation systematically covers the geometric description and kinematics of the shell, energetic principles, stress decomposition, plane stress reduction, stress resultants, and the derivation of consistent tangent moduli. Finally, the framework is extended to incorporate the effects of initial prestretch, enabling the analysis of prestrained electroelastic thin-shell structures.

\subsection{Geometric Description}
\label{subsec:geom_desc}
\begin{figure}[]
\centering
  \includegraphics[width=0.8\linewidth]{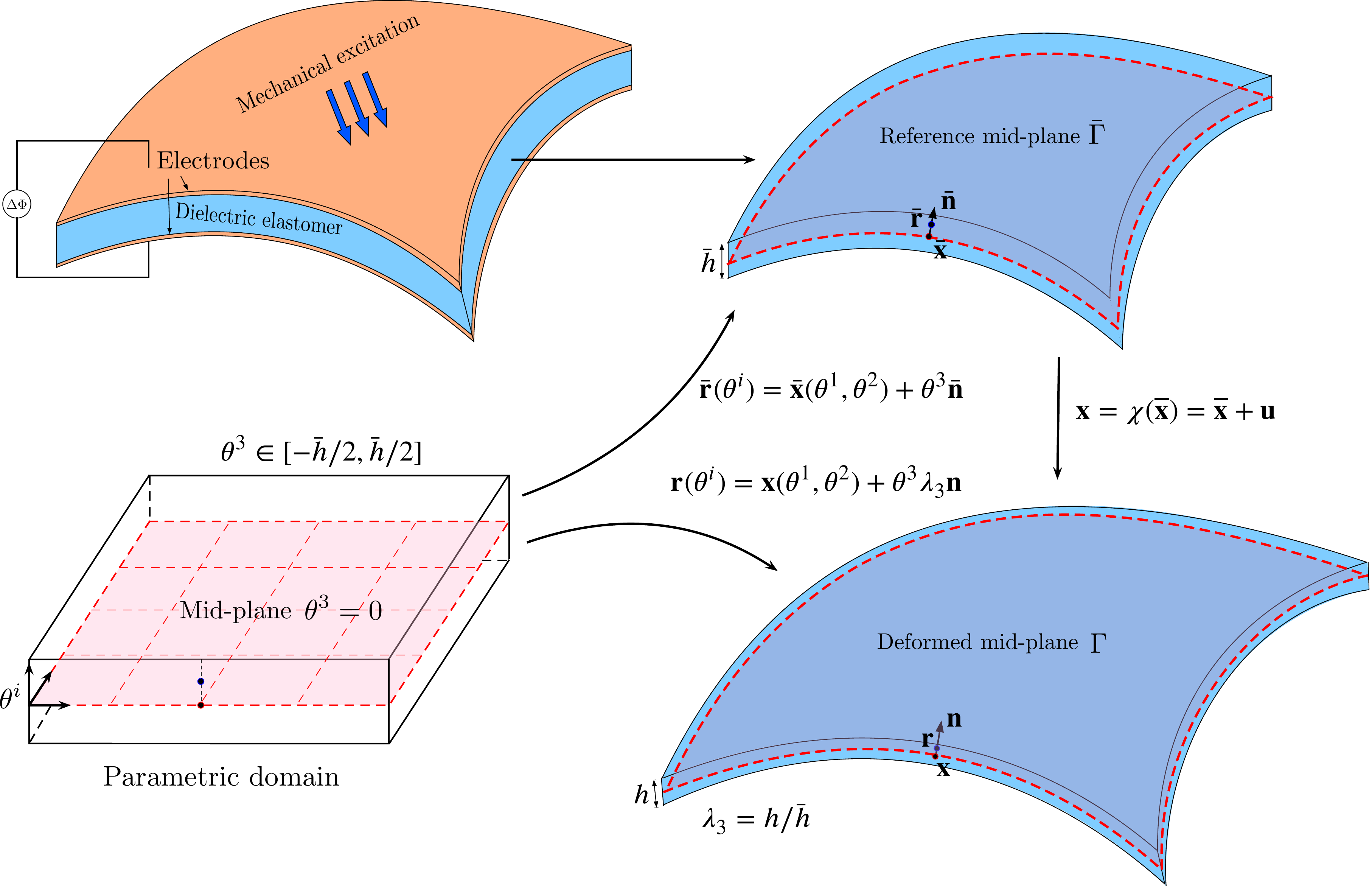}
\caption{A Kirchhoff--Love shell occupying a domain $\bar\Omega$. Each point $\bar{\mathbf{r}} \in \bar\Omega$ can be defined using quantities on the mid-surface $\bar\Gamma$ on the shell as $\bar {\mathbf r} = \bar{\mathbf x} + \theta^3 \bar{\mathbf n}$. Position vectors in the reference configuration $(\bar{\mathbf{x}} \in \bar\Gamma)$ of the mid-surface and deformed configuration $(\mathbf{x} \in \Gamma)$ of the mid-surface are related by the displacement vector $\mathbf{u}$.}
\label{fig:shell_definition}
\end{figure}
Consider a shell made of dielectric elastomer in its reference configuration as the physical domain $\bar\Omega \subset \mathbb{R}^3$, as shown in Fig.~\ref{fig:shell_definition}. The Kirchhoff--Love shell formulations is adopted to describe the mechanical behaviour of thin shell structures, which assumes that the lines that are perpendicular to the mid-surface before deformation remain straight after deformation.
Each shell point $\bar{\mathbf{r}} \in \bar\Omega$ is mapped from the parametric domain defined by the coordinates $\{ \theta^1, \theta^2, \theta^3\}$.
Assuming the shell has a uniform thickness $\bar h$ in the reference configuration, the point $\bar{\mathbf r}$ can be defined using a point on the mid-surface $\bar\Gamma$, denoted $\bar{\mathbf x} \in \bar\Gamma$, and the associated unit normal vector $\bar{\mathbf n}$ as
 \begin{equation}
\bar {\mathbf r} (\theta^1,\theta^2,\theta^3)= \bar{\mathbf x}(\theta^1,\theta^2) + \theta^3 \bar{\mathbf n}(\theta^1,\theta^2),
\label{eq: r x n relation}
 \end{equation}
 where $\theta^3 \in [-{\bar h}/{2}, {\bar h}/{2}]$.

Both the reference and deformed configurations of the shell mid-surface are mapped from the mid-surface of the parametric domain.
The corresponding mid-surface points in the reference and deformed configurations are denoted by $\bar{\mathbf x}$ and $\mathbf x$, respectively.
The position vector of a mid-surface point in the deformed configuration, $\mathbf x$, is related to its counterpart in the reference configuration, $\bar{\mathbf x}$, through
\begin{equation}
\mathbf x = \bar{\mathbf x} + \mathbf u,
\label{eq:reference_to_deformed}
\end{equation}
where $\mathbf u$ denotes the displacement vector of the mid-surface.
Moreover, the covariant basis vectors on the mid-surface in the reference and deformed configurations are defined by
\begin{equation}
\bar{\mathbf a}_\alpha = \frac{\partial \bar{\mathbf x}}{\partial \theta^\alpha}
\quad \text{and} \quad
{\mathbf a}_\alpha = \frac{\partial {\mathbf x}}{ \partial \theta^\alpha}.
\end{equation}
The corresponding unit normal vectors are then given by
\begin{equation}
\bar{\mathbf n} = \bar{\mathbf a}^3
=
\frac{\bar{\mathbf a}_1 \times \bar{\mathbf a}_2}{\bar{J}}
\quad \text{and} \quad
{\mathbf n} = {\mathbf a}^3
=
\frac{{\mathbf a}_1 \times {\mathbf a}_2}{J},
\end{equation}
where $\bar{J}$ and $J$ denote the mid-surface Jacobians in the reference and deformed configurations, respectively, defined by
\begin{equation}
\bar{J} = |\bar{\mathbf a}_1 \times \bar{\mathbf a}_2|
\quad \text{and} \quad
J = |{\mathbf a}_1 \times {\mathbf a}_2|.
\end{equation}

The covariant components of the metric tensors for the mid-surface points $\bar{\mathbf x}$ and $\mathbf x$ are respectively given by
\begin{equation}
\bar{a}_{ij} = \bar{\mathbf a}_i \cdot \bar{\mathbf a}_j \quad \text{and} \quad {a}_{ij} = {\mathbf a}_i \cdot {\mathbf a}_j.
\end{equation}
The corresponding contravariant metric tensors $\bar{a}^{ik}$ and $a^{ik}$ are defined by
\begin{equation}
\bar{a}^{ik}\bar{a}_{kj} = \delta^{i}_{j} \quad \text{and} \quad a^{ik}a_{kj} = \delta^{i}_{j},
\label{eq:co_and_contra_metric}
\end{equation}
where $\delta^{i}_{j}$ denotes the Kronecker delta. 

The thickness stretch $\lambda_3$ for a finitely deformed shell is defined by
\begin{equation}
\lambda_3 = \frac{h}{\bar h},
\end{equation}
where $h(\theta^1,\theta^2)$ is the shell thickness in the deformed configuration.
We introduce a vector $\mathbf d$ combining the thickness stretch and normal vector as
\begin{equation}
    \mathbf d = \lambda_3 \mathbf{a}_3,
\end{equation}
to write the position vector $\mathbf r$ of a point in the deformed configuration of the shell-space as
\begin{equation}
     {\mathbf r}(\theta^1,\theta^2,\theta^3) = {\mathbf x}(\theta^1,\theta^2) + \theta^3 {\mathbf d}(\theta^1,\theta^2).
\end{equation}
Thus, the three-dimensional covariant basis vectors in the shell-space of the reference and the deformed configurations, respectively, follow as
\begin{equation}
\bar{\mathbf g}_\alpha = \frac{\partial \bar{\mathbf r}}{\partial \theta^\alpha} = \bar{\mathbf a}_\alpha + \theta^3 \bar{\mathbf a}_{3,\alpha},  \quad \bar{\mathbf g}_3 = \frac{\partial \bar{\mathbf r}}{\partial \theta^3} = \bar{\mathbf a}_3,
\label{eq:cov_tensors_r}
\end{equation}
and 
\begin{equation}
{\mathbf g}_\alpha = \frac{\partial {\mathbf r}}{\partial \theta^\alpha} = {\mathbf a}_\alpha + \theta^3 \mathbf{d}_{,\alpha}, \quad 
{\mathbf g}_3 = \frac{\partial {\mathbf r}}{\partial \theta^3} = {\mathbf d}.
\label{eq:cov_tensors_d}
\end{equation}

The components of the covariant metric tensors in the shell-space are given by
\begin{equation}
\bar{g}_{ij} = \bar{\mathbf g}_i \cdot \bar{\mathbf g}_j\quad \text{and} \quad g_{ij} = \mathbf g_i \cdot \mathbf g_j,
\end{equation}
and the contravariant components of the metric tensor at point $\mathbf r$ follow as
\begin{equation}
    \bar{g}^{ij} = \bar{\mathbf g}^{i} \cdot \bar{\mathbf g}^{j} \quad \text{and} \quad {g}^{ij} = {\mathbf g}^{i} \cdot {\mathbf g}^{j},
    \label{eq:contra_metric}
\end{equation}
where $\bar{\mathbf{g}}^i$ and ${\mathbf{g}}^i$ denotes the contravariant basis vectors in reference and deformed configuration of the shell-space defined by
\begin{equation}
\bar{\mathbf{g}}^i \cdot \bar{\mathbf{g}}_j = \delta^i_j \quad \text{and} \quad {\mathbf{g}}^i \cdot {\mathbf{g}}_j = \delta^i_j.
\label{eq:deform_tensor}
\end{equation}
\subsection{Kinematics}
\label{sec:shell_kinematics}
The deformation gradient tensor $\mathbf{F}$ is defined by
\begin{equation}
    \mathbf{F} = \frac{\partial \mathbf{r}}{\partial \bar{\mathbf{r}}}= \frac{\partial\mathbf{r}}{\partial\theta^{i}}\otimes\frac{\partial\theta^{i}}{\partial\bar{\mathbf{r}}}=\mathbf{g}_{i}\otimes\mathbf{\bar{g}}^{i},
\end{equation}
thus the right Cauchy-Green deformation tensor is computed as
\begin{equation}
\mathbf{C}=\mathbf{F}^\mathrm{T}\mathbf{F}=g_{ij}\mathbf{\bar{g}}^i\otimes\mathbf{\bar{g}}^j.
\end{equation}
For thin shells undergoing moderate deformations, the out-of-plane shear terms ($C_{\alpha 3}$) are negligible. This simplifies the tensor to a membrane-dominated form:
\begin{equation}
\mathbf{C} = g_{\alpha\beta} \bar{\mathbf{g}}^\alpha \otimes \bar{\mathbf{g}}^\beta + [\lambda_3]^2 \bar{\mathbf{g}}^3 \otimes \bar{\mathbf{g}}^3.
\label{eq:right_cauchy_green_simplified}
\end{equation}
Consequently, the inverse right Cauchy-Green tensor adopts a simplified block-diagonal form:
\begin{equation}
[\mathbf{C}^{-1}] = \begin{bmatrix}
C^{-1}_{11} & C^{-1}_{12} & 0 \\
C^{-1}_{12} & C^{-1}_{22} & 0 \\
0 & 0 & C^{-1}_{33}
\end{bmatrix},
\label{eq:block_diagonal}
\end{equation}
where $C^{-1}_{33} = [\lambda_3]^{-2}$.
{
The Green-Lagrange strain tensor $\boldsymbol{\mathcal{E}}$ can be expressed as
\begin{equation}
    \boldsymbol{\mathcal{E}}=\frac{1}{2}[\mathbf{C}-\mathbf{I}]=\underbrace{\frac{1}{2}[g_{\alpha\beta}-\bar{g}_{\alpha\beta}]\bar{\mathbf{g}}^\alpha\otimes\bar{\mathbf{g}}^\beta}_{\tilde{\boldsymbol{\mathcal{E}}}}+\frac{1}{2}\left[[\lambda_3]^2-1\right]\bar{\mathbf{g}}^3\otimes\bar{\mathbf{g}}^3,
    \label{eq:green_lagrange_strain}
\end{equation}
where $\tilde{\boldsymbol{\mathcal{E}}}$ is the in-plane strain tensor. It is further decomposed into two parts:
\begin{equation}
    \tilde{\boldsymbol{\mathcal{E}}} = \boldsymbol{\epsilon} + \theta^3 \boldsymbol{\kappa},
    \label{eq:strain_decomp}
\end{equation}}
where $\boldsymbol{\epsilon}$ represents the membrane strain tensor, capturing in-plane stretching and shearing deformations, while $\boldsymbol{\kappa}$ represents the bending strain tensor, describing curvature changes due to bending or twisting. Their components are computed as
\begin{equation}
 \epsilon_{\alpha\beta} = \frac{1}{2} [a_{\alpha \beta} - \bar{a}_{\alpha \beta}]\quad \text{and} \quad \kappa_{\alpha\beta} =[- b_{\alpha \beta} + \bar {b}_{\alpha \beta}],
 \label{eq:GL_strian_2}
\end{equation}
with
\begin{equation}
\bar{a}_{\alpha\beta} = \bar{\mathbf{a}}_\alpha \cdot \bar{\mathbf{a}}_\beta,\, {a}_{\alpha\beta} = {\mathbf{a}}_\alpha \cdot {\mathbf{a}}_\beta\quad \text{and} \quad \bar{b}_{\alpha\beta} = \bar{\mathbf{a}}_{\alpha,\beta} \cdot \bar{\mathbf a}_3,\,{b}_{\alpha\beta} = {\mathbf{a}}_{\alpha,\beta} \cdot {\mathbf a}_3.
\end{equation}

\subsection{Energetic Formulation and Weak Form}
\label{subsec:weak_formulation}

The governing equations for the electromechanical equilibrium of the shell are derived from the principle of stationary potential energy. The total potential energy $\Pi_{\text{tot}}$ of the system consists of the internal energy $\Pi_{\text{int}}$ and the external work $\Pi_{\text{ext}}$:

\begin{equation}
    \Pi_{\text{tot}} = \Pi_{\text{int}} + \Pi_{\text{ext}}.
\end{equation}

{\subsubsection{Internal Energy Functional}
The internal energy accounts for the stored mechanical and electrical energy, and enforces the incompressibility constraint via a Lagrange multiplier $\tilde{p}_0$:
\begin{equation}
\Pi_{\text{int}} = \int_{\bar{\Omega}} \widetilde{W}_{\text{mech}} \, dV + \int_{\bar{\Omega}} \widetilde{W}_{\text{elec}} \, dV - \int_{\bar{\Omega}} \tilde{p}_0 [\mathcal{J} - 1] \, dV,
\label{eq:internal_energy}
\end{equation}
This functional is the direct application of the strain energy density formulation from Section~\ref{sec:strain_energy} to the principle of stationary potential energy. Equation~\eqref{eq:energy_decomposed} is adopted to separate the mechanical and electrical contributions.
$\widetilde{W}_{\mathrm{mech}}$ and $\widetilde{W}_{\mathrm{elec}}$ are the mechanical and electrical energy densities further defined in Equations~\eqref{eq:mech_energy} and \eqref{eq:elec_energy}, respectively, and $\mathcal{J} = \det(\mathbf{F}) = 1$ is the incompressibility constraint.
}
\subsubsection{First Variation and Weak Form}
The equilibrium state corresponds to a stationary point of the total energy. Taking the first variation $\delta\Pi_{\text{tot}} = 0$ yields the weak form of the balance laws:
\begin{equation}
    \delta\Pi_{\text{tot}} = \delta\Pi_{\text{int}} + \delta\Pi_{\text{ext}} = 0.
    \label{eq:weak_form}
\end{equation}
The internal virtual work is obtained as the variation of the internal energy:
\begin{equation}
    \delta\Pi_{\text{int}} = \int_{\bar{\Omega}} \delta W_{\text{int}} \, \mathrm{d}V
    = \int_{\bar{\Omega}} {\mathbf{S}} : \delta{\boldsymbol{\mathcal{E}}} \, \mathrm{d}V,
    \label{eq:variation_int_energy}
\end{equation}
where ${\mathbf{S}}$ is the total second Piola--Kirchhoff stress tensor and $\delta{\boldsymbol{\mathcal{E}}}$ is the variation of the Green-Lagrange strain tensor. 
{\textbf{Remark:} Consistent with the Kirchhoff--Love hypothesis and the plane stress assumption (${S}^{33} = 0$), the thickness strain ${\mathcal{E}}_{33}$ and transverse shear strains ${\mathcal{E}}_{\alpha 3}$ do not contribute to the virtual work. 
While the thickness stretch $\lambda_3$ is kinematically determined by the in-plane deformation and incompressibility, the vanishing stress ${S}^{33}$ ensures that ${\mathcal{E}}_{33}$ performs no work. Consequently, the internal energy depends only on the in-plane strains, which are the sole contributors to the weak form.
Thus, Equation~\eqref{eq:variation_int_energy} reduces to
\begin{equation}
    \delta\Pi_{\text{int}} = \int_{\bar{\Omega}} \tilde{\mathbf{S}} : \delta\tilde{\boldsymbol{\mathcal{E}}} \, \mathrm{d}V,
\end{equation}
where $\tilde{\mathbf{S}}$ and $\delta\tilde{\boldsymbol{\mathcal{E}}}$ are only in-plane tensors.}

The external virtual work $\delta\Pi_{\text{ext}}$ incorporates contributions from applied mechanical tractions and electrical boundary conditions, the specifics of which depend on the problem setup.

Equation~\eqref{eq:weak_form} constitutes the nonlinear variational equation to be solved. In the present isogeometric discretisation (Section~\ref{sec:discretisation}), it leads to a system of nonlinear algebraic equations for the control point displacements and electric potential.

\subsubsection{Linearisation and Material Tangent}
To solve Equation~\eqref{eq:weak_form} using the Newton-Raphson method, consistent linearisation is required. The directional derivative (linearisation) of the weak form yields the tangent stiffness operator. This involves the linearisation of the stress, which introduces the \textit{material tangent moduli} (detailed in Section~\ref{sec:material_tangent_moduli}).

The final discrete tangent stiffness matrix assembled from the finite element discretisation therefore comprises both geometric stiffness contributions (from the linearisation of the strain variation $\delta\boldsymbol{\mathcal{E}}$) and material stiffness contributions (from the tangent moduli $\mathbb{C}^{ijkl}$). The specific expressions for the condensed plane-stress tangent moduli $\hat{\mathbb{C}}^{\alpha\beta\gamma\delta}$, which are used in the shell resultant formulation, will be provided in Section~\ref{sec:material_tangent_moduli_3}.
\subsection{Total Stress in Electroelastic Thin Shells}
\label{subsec:stress}
The total stress within the dielectric elastomer shell originates from three contributions, expressed through the constitutive relationship:
\begin{equation}
    {\mathbf{S}} ={2\frac{\partial \widetilde{W}_{\mathrm{mech}}}{\partial \mathbf{C}}}+ {2\frac{\partial \widetilde{W}_{\mathrm{elec}}}{\partial \mathbf{C}}} - {\tilde{p}_0 \mathbf{C}^{-1}}.
    \label{eq:total_stress_decomposition}
\end{equation}
where the first term is the hyperelastic mechanical stress arising from the finite deformation of the dielectric elastomer. The second term corresponds to the Maxwell stress generated by the interaction between the applied electric field and dielectric material, and the third term accounts for the hydrostatic pressure that maintains the incompressibility constraint.
\subsubsection{Mechanical Stress Contribution}
The present framework is general and can be applied to any hyperelastic constitutive equation. Here, the Mooney–Rivlin model is employed as a representative example. Assuming a mechanical energy density $\widetilde{W}_{\mathrm{mech}}$ for finite deformations, the Mooney–Rivlin form is:
\begin{equation}
    \widetilde{W}_{\mathrm{mech}}(\mathbf{C}) = c_1[I_1 - 3] + c_2[I_2 - 3],
    \label{eq:mech_energy}
\end{equation}
where the invariants $I_1$ and $I_2$ are defined in Equation~\eqref{eq:mechanical_invariants}.
The corresponding stress contribution derives from the derivative of $\widetilde{W}_{\mathrm{mech}}$ with respect to $\mathbf{C}$:
\begin{equation}
    {\mathbf{S}}_{\text{mech}} \equiv 2\frac{\partial \widetilde{W}_{\mathrm{mech}}}{\partial \mathbf{C}} = 2c_1\frac{\partial I_1}{\partial \mathbf{C}} + 2c_2\frac{\partial I_2}{\partial \mathbf{C}}.
    \label{eq:mech_stress_deriv}
\end{equation}
\subsubsection{Electrically Induced Stress} 
\label{sec:electrically_induced_stress}
The electric energy density $\widetilde{W}_{\mathrm{elec}}$ captures the energy stored in the dielectric material due to polarisation under an applied electric field. For an isotropic voltage‑controlled elastomers, we recall the expression~\eqref{eq:elec_energy_1} here:
\begin{equation}
    \widetilde{W}_{\mathrm{elec}}(\mathbf{C},\bar{\mathbf{E}}) =  -\frac{1}{2}\epsilon [\bar{\mathbf{E}}\otimes \bar{\mathbf{E}}]: \mathbf{C}^{-1}.
    \label{eq:elec_energy}
\end{equation}

The electrically induced stress contribution to the total stress, the Maxwell stress, is derived as the work conjugate to the material strain measure. Applying the chain rule, its derivative the Maxwell stress yields:
\begin{equation}
    {\mathbf{S}}_{\text{elec}} \equiv 2\frac{\partial \widetilde{W}_{\mathrm{elec}}}{\partial \mathbf{C}} = \epsilon \mathbf{C}^{-1}[\bar{\mathbf{E}} \otimes \bar{\mathbf{E}}]\mathbf{C}^{-1}.
    \label{eq:elec_stress_deriv}
\end{equation}
This expression represents the Maxwell stress in material coordinates, which arises from electrostatic interactions within the dielectric medium.
For thin shell applications, a key simplification occurs when an electric potential difference $\Delta\Phi$ is applied across the thickness. The spatial electric field then simplifies to:
\begin{equation}
    \mathbf{E} = -\frac{\Delta\Phi}{h}\mathbf{a}_3 = -\frac{\Delta\Phi}{\lambda_3 \bar h}\mathbf{a}_3,
\end{equation}
The corresponding material electric field, computed via the inverse deformation gradient, becomes:
\begin{equation}
 \bar{\mathbf{E}} = \mathbf{F}^{-\mathrm{T}}\mathbf{E} = -\frac{\Delta\Phi}{\bar{h}}\bar{\mathbf{a}}_3 .
\label{eq:electric_field_shell}
\end{equation}
This formulation confirms that the electric field remains aligned with the out-of-plane direction throughout deformation.
Substituting the material electric field expression into Equation \eqref{eq:elec_stress_deriv} yields the simplified Maxwell stress tensor:
\begin{equation}
    {\mathbf{S}}_{\text{elec}} = \epsilon\left[\frac{\Delta\Phi}{\bar{h}}\right]^2 \mathbf{C}^{-1}[\bar{\mathbf{a}}_3 \otimes \bar{\mathbf{a}}_3]\mathbf{C}^{-1}.
    \label{eq:se_tensor}
\end{equation}
This expression is valid for the general anisotropic case and captures the full electromechanical coupling.

Applying the transverse isotropy simplification to the Maxwell stress expression, the quadratic product $\mathbf{C}^{-1}[\bar{\mathbf{a}}_3 \otimes \bar{\mathbf{a}}_3]\mathbf{C}^{-1}$ reduces to:
\begin{equation}
    \mathbf{C}^{-1}[\bar{\mathbf{a}}_3 \otimes \bar{\mathbf{a}}_3]\mathbf{C}^{-1} = [\lambda_3]^{-4} [\bar{\mathbf{a}}_3 \otimes \bar{\mathbf{a}}_3].
    \label{eq:out_of_plane_term}
\end{equation}
This dramatic simplification holds because  $\bar{\mathbf{a}}_3$ is an eigenvector of $\mathbf{C}^{-1}$, a property guaranteed by the assumption of transverse isotropy, satisfying $\mathbf{C}^{-1}\bar{\mathbf{a}}_3 = [\lambda_3]^{-2}\bar{\mathbf{a}}_3$. 
Exploiting the orthogonality of the shell director $\bar{\mathbf{a}}_3$ to the mid-surface, the primary stress component normal to the mid-surface is obtained by taking the tensor contraction:
\begin{equation}
    S_{\text{elec}}^{33} = [\bar{\mathbf{a}}_3 \otimes \bar{\mathbf{a}}_3] : \mathbf{S}_{\text{elec}} = \epsilon [\lambda_3]^{-4}\left[\frac{\Delta\Phi}{\bar{h}}\right]^2.
    \label{eq:s33_result}
\end{equation}
This is the principal electromechanical stress component driving thickness changes in actuation.
The transverse isotropy assumption and normal electric field orientation cause all off-diagonal and in-plane components of $\mathbf{S}_{\text{elec}}$ to vanish. Mathematically, this occurs because:
\[
{S}_{\text{elec}}^{\alpha\beta} = [\bar{\mathbf{a}}_\alpha \otimes \bar{\mathbf{a}}_\beta] : \mathbf{S}_{\text{elec}} = 0 \quad \text{for} \quad \alpha,\beta = 1,2,
\]
due to orthogonality between $\bar{\mathbf{a}}_3$ and $\bar{\mathbf{a}}_\alpha$. Physically, this suppression of in-plane electromechanical coupling arises because the applied electric field is oriented exclusively normal to the mid-surface.

\subsubsection{Explicit Enforcement of Plane Stress and Incompressibility}
For thin-shell structures, three physical considerations justify the plane stress assumption:
\begin{enumerate}
    \item \textit{Dimensional disparity}: The thickness dimension is orders of magnitude smaller than in-plane dimensions
    \item \textit{Boundary conditions}: Both top and bottom surfaces are traction-free ($\mathbf{t} = \mathbf{0}$)
    \item \textit{Stress magnitude}: Through-thickness stresses are negligible compared to in-plane stresses
\end{enumerate}
{
\textbf{Remark:} The Maxwell stress is an internal electromechanical coupling effect and does not constitute an external mechanical traction on the boundaries. The electrodes impose only the electric potential, leaving the mechanical traction $\mathbf{t}=\mathbf{0}$. The thickness component of the Maxwell stress is internally balanced by the Lagrange multiplier $\tilde{p}_0$ to satisfy ${S}^{33}=0$.
}

These conditions remain valid for electroelastic shells under actuation, leading to the equilibrium condition:
\begin{equation}
    {S}^{33} = S^{33}_{\text{elec}} + S^{33}_{\text{mech}} - \tilde{p}_0 C^{33} = 0,
    \label{eq:plane_stress_condition}
\end{equation}
where $C^{33} = [\lambda_3]^{-2}$. This equation expresses the vanishing normal stress in the thickness direction.
The Lagrange multiplier $\tilde{p}_0$ is explicitly determined by solving Equation \eqref{eq:plane_stress_condition}. The resulting expression is:
\begin{align}
    \tilde{p}_0 &= [\lambda_3]^2\left[S^{33}_{\text{mech}} + S^{33}_{\text{elec}}\right] \nonumber \\
    &= 2[\lambda_3]^2\frac{\partial \widetilde{W}_{\mathrm{mech}}}{\partial C_{33}} + [\lambda_3]^{-2} \epsilon \left[\frac{\Delta\Phi}{\bar{h}}\right]^2.
    \label{eq:pressure_expression}
\end{align}
This solution strategy eliminates  $\tilde{p}_0$ as an additional unknown variable while ensuring the exact enforcement of both constraints.
It clearly separates into mechanical and electrical contributions here.
The electrical term $[\lambda_3]^{-2} \epsilon [\Delta\Phi/\bar{h}]^2$ represents `electrostatic pressure', which decreases with thickness stretch.

Substituting $\tilde{p}_0$ into the general stress expression yields the working in-plane stress components:
\begin{align}
    \tilde{S}^{\alpha\beta} 
    &= \underbrace{2\left[\frac{\partial \widetilde{W}_{\mathrm{mech}}}{\partial C_{\alpha\beta}} - [\lambda_3]^2\frac{\partial \widetilde{W}_{\mathrm{mech}}}{\partial C_{33}}C^{\alpha\beta}\right]}_{\text{Mechanical stress}} - \underbrace{[\lambda_3]^{-2} \epsilon \left[\frac{\Delta\Phi}{\bar{h}}\right]^2C^{\alpha\beta}}_{\text{Electrically induced stress}}
    \label{eq:inplane_stress}
\end{align}
This formulation reveals that the electrical term scales with $(\Delta\Phi)^2$ and $C^{\alpha\beta}$, creating voltage-dependent in-plane stresses and they increase with thickness reduction. The decoupled structure enables efficient implementation while capturing the essential electromechanical coupling mechanisms.

\subsubsection{Stress Resultants for Electroelastic Thin Shells}
\label{sec:Resultants}
The thin shell formulation reduces the three-dimensional continuum to a two-dimensional surface with a thickness. The internal forces are expressed as the resultants of the integrated stress through the thickness, defined by the in-plane second Piola--Kirchhoff stress $\tilde{\mathbf{S}}$. The resultant of the in-plane stress $\boldsymbol{\mathcal{N}}$ and the bending moment $\boldsymbol{\mathcal{M}}$ have their components calculated as:
\begin{subequations}
\begin{align}
n^{\alpha\beta} &= \int_{-\frac{\bar{h}}{2}}^{\frac{\bar{h}}{2}} \tilde{S}^{\alpha\beta} J_c \, \mathrm{d}\theta^3, \\
m^{\alpha\beta} &= \int_{-\frac{\bar{h}}{2}}^{\frac{\bar{h}}{2}} \tilde{S}^{\alpha\beta} \theta^3 J_c \, \mathrm{d}\theta^3,
\end{align}
\label{eq:stress_resultants}
\end{subequations}
where $\theta^3 \in [-\bar{h}/2, \bar{h}/2]$ is the coordinate along the reference thickness and $J_c$ is the thickness-direction Jacobian correction factor:
\begin{equation}
J_c = \frac{\left| [\bar{\mathbf{g}}_1 \times \bar{\mathbf{g}}_2] \cdot \bar{\mathbf{g}}_3 \right|}{\left| [\bar{\mathbf{a}}_1 \times \bar{\mathbf{a}}_2] \cdot \bar{\mathbf{a}}_3 \right|},
\label{eq:jacobian_correction}
\end{equation}
accounting for volume changes between the reference base vectors $\bar{\mathbf{g}}_i$ and mid-surface basis $\bar{\mathbf{a}}_\alpha$.

The stress increments are linearised via the constitutive relation between ${\mathbf{S}}$ and the Green-Lagrange strain $\boldsymbol{\mathcal{E}}$ (Equation~\eqref{eq:total_stress_decomposition}) as
\begin{equation}
\mathrm{d} {S}^{ij}=\frac{\partial S^{ij}}{\partial \mathcal{E}_{kl}}\mathrm{d}\mathcal{E}_{kl}=2\frac{\partial S^{ij}}{\partial C_{kl}}\mathrm{d}\mathcal{E}_{kl}=\hat{\mathbb{C}}^{ijkl}\mathrm{d}\mathcal{E}_{kl}
\label{eq:stress_increment}
\end{equation}
where $\hat{\mathbb{C}}^{ijkl}$ are the components of the fourth-order elasticity tensor. 
\subsection{Variation of Total Energy}
The total potential energy $\Pi_{\text{tot}}$ of the electroelastic thin shell system comprises internal energy from finite deformation and polarisation ($\Pi_{\text{int}}$) and external work contributions ($\Pi_{\text{ext}}$). Its first variation is
\begin{equation}
\delta\Pi_{\text{tot}} = \delta\Pi_{\text{int}} + \delta\Pi_{\text{ext}} = \int_{\bar{\Omega}} \delta W_{\text{int}} \mathrm{~d} V + \int_{\Omega}\delta  W_{\text{ext}} \mathrm{~d} V.
\end{equation}
The tangent stiffness required for Newton-Raphson iterations derives from consistent linearisation of this variation, are detailed in Section~\ref{sec:weak_form_discretisation}.

\subsection{Material Tangent Moduli and Static Condensation}
\label{sec:material_tangent_moduli}
\subsubsection{Pressure Derivatives for Constraint Enforcement}
The Lagrange multiplier $\tilde{p}_0$, which enforces incompressibility ($\mathcal{J} = 1$), depends implicitly on deformation through the plane stress condition~\eqref{eq:pressure_expression}. Its derivatives with respect to the deformation components are essential for consistent tangent moduli:
\begin{align}
    \frac{\partial \tilde{p}_0}{\partial C_{\alpha\beta}} &= 2[\lambda_3]^2\frac{\partial^2 \widetilde{W}_{\mathrm{mech}}}{\partial C_{33}\partial C_{\alpha\beta}}, \\
    \frac{\partial \tilde{p}_0}{\partial C_{33}} &= 2[\lambda_3]^2\frac{\partial^2 \widetilde{W}_{\mathrm{mech}}}{\partial C_{33}^2} + 2\frac{\partial \widetilde{W}_{\mathrm{mech}}}{\partial C_{33}} - [\lambda_3]^{-4} \epsilon \left[\frac{\Delta\Phi}{\bar{h}}\right]^2.
\end{align}
These capture how pressure responds to:  
(1) in-plane stretching ($\alpha,\beta=1,2$) through mechanical-kinematic coupling, and  
(2) thickness changes ($C_{33}$) with explicit electromechanical contributions.

\subsubsection{Material Tangent Moduli Derivation}
The fourth-order elasticity tensor $\mathbb{C}^{ijkl} \equiv 2\frac{\partial {S}^{ij}}{\partial C_{kl}}$ is derived by consistent differentiation of the total stress (Equation~\eqref{eq:total_stress_decomposition}) with respect to the right Cauchy-Green tensor. Its partitions exhibit distinct symmetries:

\paragraph*{In-plane Moduli ($\alpha\beta\gamma\delta$)}
\begin{align}
    \mathbb{C}^{\alpha\beta\gamma\delta} &= 4\pderiv{\widetilde{W}_{\mathrm{mech}}}{C_{\alpha\beta}}{C_{\gamma\delta}} - 4[\lambda_3]^2\left[\pderiv{\widetilde{W}_{\mathrm{mech}}}{C_{33}}{C_{\gamma\delta}}C^{\alpha\beta} + \pderiv{\widetilde{W}_{\mathrm{mech}}}{C_{33}}{C_{\alpha\beta}}C^{\gamma\delta}\right] \nonumber \\
    &\quad - {\left[2[\lambda_3]^2\pder{\widetilde{W}_{\mathrm{mech}}}{C_{33}} + [\lambda_3]^{-2} \epsilon \left[\frac{\Delta\Phi}{\bar{h}}\right]^2\right]}\left[C^{\alpha\beta}C^{\gamma\delta} - C^{\alpha\gamma}C^{\beta\delta} - C^{\alpha\delta}C^{\beta\gamma}\right],
    \label{eq:C_alpha_beta_gamma_delta}
\end{align}
where the last term arises from the product rule applied to $-\tilde{p}_0\mathbf{C}^{-1}$. This expression has the major symmetry $\mathbb{C}^{\alpha\beta\gamma\delta} = \mathbb{C}^{\gamma\delta\alpha\beta}$.

\paragraph*{Thickness Coupling Moduli ($\alpha\beta33$)}
\begin{align}
    \mathbb{C}^{\alpha\beta33} &= -C^{\alpha\beta}\left[6\pder{\widetilde{W}_{\mathrm{mech}}}{C_{33}} + 4[\lambda_3]^2\pdersq{\widetilde{W}_{\mathrm{mech}}}{C_{33}} - [\lambda_3]^{-4} \epsilon \left[\frac{\Delta\Phi}{\bar{h}}\right]^2\right],
    \label{eq:C_alpha_beta_33}
\end{align}
quantifying how in-plane stresses change with thickness stretch.

\paragraph*{Thickness Moduli ($3333$)}
\begin{align}
    \mathbb{C}^{3333} &= -[\lambda_3]^{-2}\left[6\pder{\widetilde{W}_{\mathrm{mech}}}{C_{33}} + 4[\lambda_3]^2\pdersq{\widetilde{W}_{\mathrm{mech}}}{C_{33}} - [\lambda_3]^{-4} \epsilon \left[\frac{\Delta\Phi}{\bar{h}}\right]^2\right],
    \label{eq:C_3333}
\end{align}
governing thickness-direction stiffness. Its negative definiteness reflects the kinematic constraint.

\subsubsection{Static Condensation for Plane Stress}
\label{sec:material_tangent_moduli_3}
To enforce the plane stress condition ($\tilde{S}^{33}=0$), the moduli undergo static condensation:
\begin{align}
    \hat{\mathbb{C}}^{\alpha\beta\gamma\delta} &= \mathbb{C}^{\alpha\beta\gamma\delta} - \frac{\mathbb{C}^{\alpha\beta33}\mathbb{C}^{33\gamma\delta}}{\mathbb{C}^{3333}}.
    \label{eq:condensation}
\end{align}
This simplifies the 3D constitutive relation to 2D, eliminating explicit $\tilde{p}_0$ dependence. Substituting Eqs. \eqref{eq:C_alpha_beta_gamma_delta}-\eqref{eq:C_3333} yields:
\begin{align}
    \hat{\mathbb{C}}^{\alpha\beta\gamma\delta} &= 4\pderiv{\widetilde{W}_{\mathrm{mech}}}{C_{\alpha\beta}}{C_{\gamma\delta}} - 4[\lambda_3]^2\left[\pderiv{\widetilde{W}_{\mathrm{mech}}}{C_{33}}{C_{\gamma\delta}}C^{\alpha\beta} + \pderiv{\widetilde{W}_{\mathrm{mech}}}{C_{33}}{C_{\alpha\beta}}C^{\gamma\delta}\right] \nonumber \\
    &\quad + C^{\alpha\beta}C^{\gamma\delta}\left[6[\lambda_3]^2\pder{\widetilde{W}_{\mathrm{mech}}}{C_{33}} + 4[\lambda_3]^4\pdersq{\widetilde{W}_{\mathrm{mech}}}{C_{33}} - [\lambda_3]^{-2} \epsilon \left[\frac{\Delta\Phi}{\bar{h}}\right]^2\right],
    \label{eq:reduced_tangent}
\end{align}
where $\hat{\mathbb{C}}^{\alpha\beta\gamma\delta}$ is the plane-stress-reduced tangent modulus, which incorporates both mechanical and electrical effects while satisfying plane stress constraints intrinsically. 
\subsection{Prestretch in Electroelastic Thin Shells}
\label{sec:prestrain_theory}
\begin{figure}
    \centering
    \includegraphics[width=\linewidth]{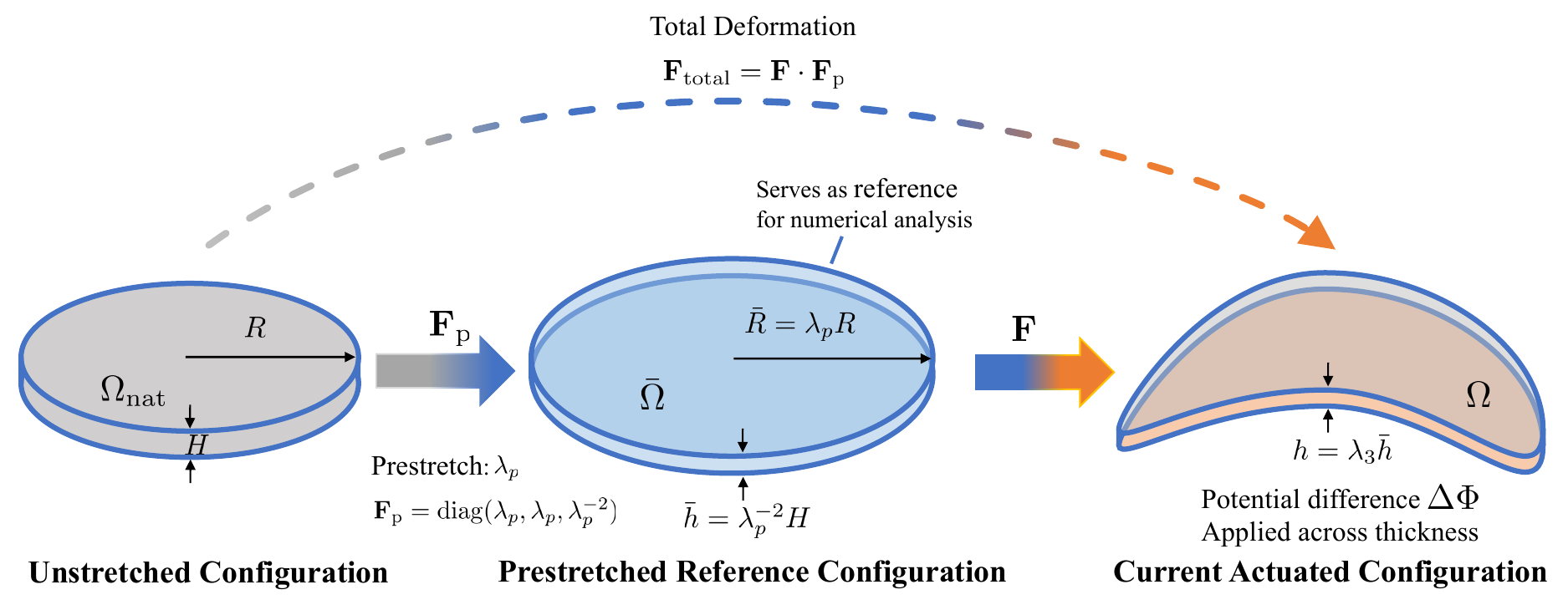}
    \caption{Progressive configurations for prestretched electroelastic thin shells}
    \label{fig:prestretched_plate}
\end{figure}
In manufacturing and application, dielectric elastomers are often subjected to prestretch to enhance actuation performance or achieve specific configurations. Figure~\ref{fig:prestretched_plate} illustrates three states for prestretched electroelastic shells. The prestretched state, denoted as $\bar\Omega$, serves as the reference configuration for subsequent electromechanical analysis. The deformation from the stress-free natural state $\Omega_{\text{nat}}$ to the prestretched reference configuration $\bar\Omega$ is described by the prestretch gradient $\mathbf{F}_\mathrm{p}$. The total deformation gradient from $\Omega_{\text{nat}}$ to the current configuration $\Omega$ decomposes multiplicatively as:
\begin{equation}
    \mathbf{F}_{\text{total}} = \mathbf{F} \cdot \mathbf{F}_\mathrm{p},
    \label{eq:F_total_decomp}
\end{equation}
where $\mathbf{F}$ is the deformation gradient from the prestretched reference configuration to the current configuration, consistent with Section~\ref{sec:shell_kinematics}.

The total right Cauchy-Green deformation tensor relative to the natural state is:

\begin{equation}
    \mathbf{C}_{\text{total}} = \mathbf{F}_\mathrm{p}^\top \mathbf{C} \mathbf{F}_\mathrm{p}.
    \label{eq:C_total_relation}
\end{equation}

The mechanical strain energy density $\widetilde{W}_{\mathrm{mech}}$ remains defined with respect to the natural state, making it a function of $\mathbf{C}_{\text{total}}$.

\subsubsection{Isotropic Prestretch}
\label{subsec:isotropic_prestretch}

A common prestretch pattern is isotropic in-plane stretching accompanied by thickness reduction to satisfy incompressibility. This is represented by
\begin{equation}
    \mathbf{F}_\mathrm{p} = \operatorname{diag}(\lambda_p, \lambda_p, \lambda_p^{-2}),
    \label{eq:Fp_isotropic}
\end{equation}
where $\lambda_p$ is the in-plane prestretch ratio. The thickness component $\lambda_p^{-2}$ ensures $\det(\mathbf{F}_\mathrm{p}) = 1$.

Under this assumption, the components of $\mathbf{C}_{\text{total}}$ relate to those of $\mathbf{C}$ (defined in the prestretched configuration) as:

\begin{equation}
    \mathbf{C}_{\text{total}} = 
    \begin{bmatrix}
        \lambda_p^2 C_{11} & \lambda_p^2 C_{12} & 0 \\
        \lambda_p^2 C_{12} & \lambda_p^2 C_{22} & 0 \\
        0 & 0 & \lambda_p^{-4} C_{33}
    \end{bmatrix},
    \label{eq:C_total_isotropic}
\end{equation}
where $C_{33} = [\lambda_3]^2$ and $\lambda_3$ is the through-thickness stretch from the prestretched reference state.

\subsubsection{Stress Resultants with Prestretch}
\label{subsec:stress_with_prestretch}

The in-plane stress components in the prestretched reference configuration follow the same general form as Equation~\eqref{eq:inplane_stress}, but with the derivatives of $\widetilde{W}_{\mathrm{mech}}$ now taken with respect to $\mathbf{C}_{\text{total}}$. For the Mooney–Rivlin model, the derivatives become:

\begin{align}
    \frac{\partial\widetilde{W}_{\mathrm{mech}}}{\partial C_{\alpha\beta}} &= 
    c_1\lambda_p^2\bar{g}^{\alpha\beta} + c_2\Bigl[\lambda_p^4\bigl[C_{\gamma\delta}\bar{g}^{\gamma\delta}\bar{g}^{\alpha\beta} - \bar{g}^{\alpha\gamma}C_{\gamma\delta}\bar{g}^{\delta\beta}\bigr] 
    + \lambda_p^{-2}\bar{g}^{\alpha\beta}C_{33}\Bigr], \label{eq:dW_dC_ab_pre} \\
    \frac{\partial\widetilde{W}_{\mathrm{mech}}}{\partial C_{33}} &= 
    c_1\lambda_p^{-4} + c_2\lambda_p^{-2}C_{\gamma\delta}\bar{g}^{\gamma\delta}. \label{eq:dW_dC_33_pre}
\end{align}

Substituting these into Equation~\eqref{eq:inplane_stress} yields the in-plane stress components explicitly accounting for prestretch:

\begin{align}
    \tilde{S}^{\alpha\beta} = 2\Bigg[ &c_1\lambda_p^2\bar{g}^{\alpha\beta} + c_2\Bigl[\lambda_p^4\bigl[C_{\gamma\delta}\bar{g}^{\gamma\delta}\bar{g}^{\alpha\beta} - \bar{g}^{\alpha\gamma}C_{\gamma\delta}\bar{g}^{\delta\beta}\bigr] + \lambda_p^{-2}\bar{g}^{\alpha\beta}C_{33}\Bigr] \nonumber \\
    &- [\lambda_3]^2 \bigl[c_1\lambda_p^{-4} + c_2\lambda_p^{-2}C_{\gamma\delta}\bar{g}^{\gamma\delta}\bigr] C^{\alpha\beta} \Bigg]
    - [\lambda_3]^{-2} \epsilon \left[ \frac{\Delta\Phi}{\bar{h}} \right]^2 C^{\alpha\beta}. \label{eq:inplane_stress_with_prestretch}
\end{align}

\subsubsection{Tangent Moduli with Prestretch}
\label{subsec:tangent_moduli_prestretch}

The consistent tangent moduli for the prestressed configuration are derived in the same manner as in Equation~\eqref{eq:reduced_tangent}, but with the derivatives of $\widetilde{W}_{\mathrm{mech}}$ now evaluated with respect to $\mathbf{C}_{\text{total}}$ and then transformed to the prestressed reference configuration. The reduced in-plane tangent moduli $\hat{\mathbb{C}}^{\alpha\beta\gamma\delta}$ are given by the same expression as Equation~\eqref{eq:reduced_tangent}, but with the following derivatives for the Mooney–Rivlin model under isotropic prestretch:
\begin{align}
    \frac{\partial^2\widetilde{W}_{\mathrm{mech}}}{\partial C_{\alpha\beta}\partial C_{\gamma\delta}} &= c_2\lambda_p^4\left[\bar{g}^{\alpha\beta}\bar{g}^{\gamma\delta} \frac{1}{2}\left[\bar{g}^{\alpha\gamma}\bar{g}^{\beta\delta} + \bar{g}^{\alpha\delta}\bar{g}^{\beta\gamma}\right]\right], \label{eq:d2W_dCab_dCcd_pre} \\
    \frac{\partial^2\widetilde{W}_{\mathrm{mech}}}{\partial C_{\alpha\beta}\partial C_{33}} &= c_2\lambda_p^{-2}\bar{g}^{\alpha\beta}, \label{eq:d2W_dCab_dC33_pre} \\
    \frac{\partial^2\widetilde{W}_{\mathrm{mech}}}{\partial C_{33}^2} &= 0. \label{eq:d2W_dC33_dC33_pre}
\end{align}

These expressions, along with the first derivatives given in Equations~\eqref{eq:dW_dC_ab_pre} and \eqref{eq:dW_dC_33_pre}, are used to compute the tangent moduli in the prestressed configuration. The electrical contribution remains unchanged, as it is independent of the prestretch.

\subsubsection{Remarks on Implementation}
\label{subsec:prestrain_implementation_remarks}

Incorporating prestretch within the isogeometric shell framework requires the following modifications:
\begin{enumerate}
    \item The reference geometry is defined in the prestretched configuration $\bar\Omega$.
    \item The prestretch tensor $\mathbf{F}_\mathrm{p}$ is stored as a field (constant or varying spatially).
    \item The strain energy derivatives in the weak form and tangent stiffness are evaluated using Equations~\eqref{eq:dW_dC_ab_pre}--\eqref{eq:d2W_dC33_dC33_pre}.
    \item The thickness stretch $\lambda_3$ is determined from the incompressibility constraint $\mathcal{J}=1$, which now includes the prestretch effect.
\end{enumerate}
This formulation enables the analysis of prestrained dielectric elastomer shells undergoing large electromechanical deformations while maintaining the $C^1$-continuity requirements of the Kirchhoff–Love shell theory.
\section{Numerical Implementation}
\label{sec:implementation}
The nonlinear finite element implementation leverages the \textsc{deal.II} library~\cite{Arndt2022-dealii94} to realize the electroelastic thin shell formulation.
The primary numerical challenge lies in the coupling between the large-deformation mechanics of a Kirchhoff–Love shell and the electrostatic forces arising from an applied voltage. 
This section outlines the key components of the numerical discretisation, the linearisation of the weak form, and the resulting solution procedure.

\subsection{Subdivision Surface Discretisation}
\label{sec:discretisation}
As established in our previous work on hyperelastic thin shells \cite{liu2024computational}, the mid-surface geometry and the displacement field are discretised using Catmull–Clark subdivision surfaces. This choice is motivated by the $C^1$-continuity requirement of the Kirchhoff–Love theory, which necessitates basis functions for the displacement in the Sobolev space $H^2(\Omega)$.
Subdivision surfaces provide smooth, $C^1$-continuous limit surfaces everywhere, even on unstructured control meshes containing extraordinary vertices.
 The mid-surface in the reference configuration $\bar{\mathbf{x}}$ and the displacement field $\mathbf{u}$ are approximated using the same set of subdivision basis functions
\begin{align}
\bar{\mathbf{x}}(\theta^1, \theta^2) &\approx \sum_{A=1}^{n_{\text{node}}} N^A(\theta^1, \theta^2) \bar{\mathbf{X}}^A, \\
\mathbf{u}(\theta^1, \theta^2) &\approx \sum_{A=1}^{n_{\text{node}}} N^A(\theta^1, \theta^2) \mathbf{u}^A,
\label{eq:displacement_interpolation}
\end{align}
where $\bar{\mathbf{X}}^A$ and $\mathbf{u}^A$ are the reference position and displacement vector of control point $A$, respectively, and $n_{\text{node}}$ is the number of control points in the support of a given parametric location. 
The deformed mid-surface position is then ${\mathbf{x}} = \bar{\mathbf{x}} + {\mathbf{u}}$.
The functions $N^A(\theta^1, \theta^2)$ are the cubic Catmull-Clark subdivision bases.
For regular patches (interior elements with valence 4), they are equivalent to bi-cubic B-splines. 
For patches containing extraordinary vertices, Stam's algorithm \cite{stam1998exact} is employed for fast and exact evaluation.
\subsection{Discretised Kinematic Fields}
All strain and curvature measures are derived from the discretised mid-surface geometry. The covariant basis vectors on the deformed mid-surface are computed as the partial derivatives of the interpolated position:
\begin{equation}
\mathbf{a}_\alpha = \mathbf{x}_{,\alpha} = \bar{\mathbf{a}}_\alpha + \sum_{A} N^A_{,\alpha} \mathbf{u}^A.
\end{equation}
The metric components $\epsilon_{\alpha\beta}$ and curvature components $\kappa_{\alpha\beta}$ (and their referential counterparts) are then calculated according to their definitions in Section 2.2. The Green–Lagrange membrane and bending strain components, as per Equation~\eqref{eq:GL_strian_2}, become functions of the nodal displacements:
\begin{align}
\epsilon_{\alpha\beta} &= \frac{1}{2}[\mathbf{a}_\alpha \cdot \mathbf{a}_\beta - \bar{\mathbf{a}}_\alpha \cdot \bar{\mathbf{a}}_\beta], \\
\kappa_{\alpha\beta} &= -[\mathbf{a}_{\alpha,\beta} \cdot \mathbf{a}_3]+ [\bar{\mathbf{a}}_{\alpha,\beta} \cdot \bar{\mathbf{a}}_3].
\end{align}
The thickness stretch $\lambda_3$, required for the stress evaluation, is computed from the area change of the in-plane basis vectors:
\begin{equation}
\lambda_3 = \frac{|\bar{\mathbf{a}}_1 \times \bar{\mathbf{a}}_2|}{|{\mathbf{a}}_1 \times {\mathbf{a}}_2|} = \frac{\bar{J}}{J}.
\end{equation}
The variations of the strains, $\delta\epsilon_{\alpha\beta} $ and $\delta\kappa_{\alpha\beta}$, are obtained by taking the directional derivative (or Gateaux derivative) of the above expressions with respect to the displacement field $\mathbf{u}$, leading to expressions linear in $\delta \mathbf{u}$.
\subsection{Weak Form and Consistent Linearisation}
\label{sec:weak_form_discretisation}
The principle of virtual work for the shell, including the internal electromechanical stresses and external loads, is
\begin{equation}
\delta \Pi_{\text{tot}} = \delta\Pi_{\text{int}} + \delta\Pi_{\text{ext}} = 0,
\end{equation}
where the internal virtual work is the integral of the stress resultants (Equation~\eqref{eq:stress_resultants}) working through the virtual strains:
\begin{equation}
\delta\Pi_{\text{int}} = \int_{\bar{\Gamma}} \left[ n^{\alpha\beta} \, \delta\epsilon_{\alpha\beta} + m^{\alpha\beta} \, \delta\kappa_{\alpha\beta} \right] d\bar{\Gamma}.
\label{eq:internal_virtual_work}
\end{equation}
Here, $\bar{\Gamma}$ is the shell mid-surface in the reference configuration, and $d\bar{\Gamma} = \bar{J} \, d\theta^1 d\theta^2$.
The stress resultants $n^{\alpha\beta}$ and $m^{\alpha\beta}$ encapsulate the full electromechanical coupling. They are computed by integrating the Piola--Kirchhoff stress $\tilde{S}^{\alpha\beta}$ (Equation~\eqref{eq:inplane_stress}) through the thickness, as defined in Equation~\eqref{eq:stress_resultants}. The integration is performed numerically using a Gauss quadrature rule along $\theta^3$. The Jacobian correction factor $J_c$ in Equation~\eqref{eq:jacobian_correction} accounts for the curvature of the shell space. The external virtual work $\delta\Pi_{\text{ext}}$ accounts for mechanical surface tractions, pressure loads, and follower forces, which are standard in shell formulations.

The solution is obtained via the Newton--Raphson method, which requires the linearisation of Equation~\eqref{eq:internal_virtual_work}. The linearised increment of the internal virtual work is:
\begin{equation}
\varDelta(\delta\Pi_{\text{int}}) = \int_{\bar{\Gamma}} \left[ \varDelta n^{\alpha\beta} \, \delta\epsilon_{\alpha\beta} + n^{\alpha\beta} \, \varDelta(\delta\epsilon_{\alpha\beta}) + \varDelta m^{\alpha\beta} \, \delta\kappa_{\alpha\beta} + m^{\alpha\beta} \, \varDelta(\delta\kappa_{\alpha\beta}) \right] d\bar{\Gamma}.
\end{equation}
The stress resultant increments $\varDelta n^{\alpha\beta}$ and $\varDelta m^{\alpha\beta}$ are related to the increments of the strain measures via the constitutive tangent moduli:
\begin{align}
\varDelta n^{\alpha\beta} &= \mathbb{A}^{\alpha\beta\gamma\delta} \, \varDelta\epsilon_{\gamma\delta} + \mathbb{B}^{\alpha\beta\gamma\delta} \, \varDelta\kappa_{\gamma\delta}, \\
\varDelta m^{\alpha\beta} &= \mathbb{B}^{\gamma\delta\alpha\beta} \, \varDelta\epsilon_{\gamma\delta} + \mathbb{D}^{\alpha\beta\gamma\delta} \, \varDelta\kappa_{\gamma\delta}.
\end{align}
The tangent stiffness tensors $\mathbb{A}$, $\mathbb{B}$, and $\mathbb{D}$ are the integrated counterparts of the reduced material tangent modulus $\hat{\mathbb{C}}^{\alpha\beta\gamma\delta}$ (Equation~\eqref{eq:reduced_tangent}):
\begin{align}
\mathbb{A}^{\alpha\beta\gamma\delta} &= \int_{-\frac{\bar{h}}{2}}^{\frac{\bar{h}}{2}} \hat{\mathbb{C}}^{\alpha\beta\gamma\delta} \, J_c \, d\theta^3, \\
\mathbb{B}^{\alpha\beta\gamma\delta} &= \int_{-\frac{\bar{h}}{2}}^{\frac{\bar{h}}{2}} \hat{\mathbb{C}}^{\alpha\beta\gamma\delta} \, \theta^3 \, J_c \, d\theta^3, \\
\mathbb{D}^{\alpha\beta\gamma\delta} &= \int_{-\frac{\bar{h}}{2}}^{\frac{\bar{h}}{2}} \hat{\mathbb{C}}^{\alpha\beta\gamma\delta} \, [\theta^3]^2 \, J_c \, d\theta^3.
\label{eq:tangent_tensors}
\end{align}
Note that $\hat{\mathbb{C}}^{\alpha\beta\gamma\delta}$ depends on the current deformation (via $\mathbf{C}$ and $\lambda_3$) and the applied voltage $\Delta\Phi$. Its explicit form, provided in Equation~\eqref{eq:reduced_tangent}, is evaluated at each integration point during the assembly of the stiffness matrix.
The terms $\varDelta(\delta\epsilon_{\alpha\beta})$ and $\varDelta(\delta\kappa_{\alpha\beta})$ arise from the geometric non-linearity and contribute to the \emph{geometric stiffness matrix}. Their detailed derivation involves the second variation of the kinematic quantities, which is standard in non-linear shell formulations (see, e.g., \cite{kiendl2015isogeometric} for the Kirchhoff-Love case and \cite{simo1990stress} for the general geometrically nonlinear shell framework).
\subsection{Finite Element Formulation}
\label{sec:FE_form}

Substituting the discretised kinematic fields into the weak form yields the discrete nonlinear equilibrium equations. For a control point $A$, the residual force is defined by
\begin{equation}
\ary{R}^{A}
=
\ary{F}^{A,\mathrm{int}}
-
\ary{F}^{A,\mathrm{ext}},
\end{equation}
where the internal force vector is given by
\begin{equation}
\ary{F}^{A,\mathrm{int}}
=
\int_{\bar{\Gamma}}
\left[
n^{\alpha\beta}
\frac{\partial \epsilon_{\alpha\beta}}
{\partial \ary{u}^{A}}
+
m^{\alpha\beta}
\frac{\partial \kappa_{\alpha\beta}}
{\partial \ary{u}^{A}}
\right]
d\bar{\Gamma}.
\label{eq:internal_force}
\end{equation}

The tangent stiffness matrix is obtained by consistent linearisation of the residual vector with respect to the displacement degrees of freedom,
\begin{equation}
\ary{K}^{AB}
=
\frac{\partial \ary{R}^{A}}
{\partial \ary{u}^{B}}.
\end{equation}

Substituting the constitutive relations for the stress resultant increments yields
\begin{align}
\ary{K}^{AB}
=
\int_{\bar{\Gamma}}
\Bigg[
&
\frac{\partial \epsilon_{\alpha\beta}}
{\partial \ary{u}^{A}}
\,
\mathbb{A}^{\alpha\beta\gamma\delta}
\,
\frac{\partial \epsilon_{\gamma\delta}}
{\partial \ary{u}^{B}}
+
\frac{\partial \epsilon_{\alpha\beta}}
{\partial \ary{u}^{A}}
\,
\mathbb{B}^{\alpha\beta\gamma\delta}
\,
\frac{\partial \kappa_{\gamma\delta}}
{\partial \ary{u}^{B}}
\nonumber\\
&
+
\frac{\partial \kappa_{\alpha\beta}}
{\partial \ary{u}^{A}}
\,
\mathbb{B}^{\gamma\delta\alpha\beta}
\,
\frac{\partial \epsilon_{\gamma\delta}}
{\partial \ary{u}^{B}}
+
\frac{\partial \kappa_{\alpha\beta}}
{\partial \ary{u}^{A}}
\,
\mathbb{D}^{\alpha\beta\gamma\delta}
\,
\frac{\partial \kappa_{\gamma\delta}}
{\partial \ary{u}^{B}}
\nonumber\\
&
+
n^{\alpha\beta}
\frac{\partial^{2}\epsilon_{\alpha\beta}}
{\partial \ary{u}^{A}\partial \ary{u}^{B}}
+
m^{\alpha\beta}
\frac{\partial^{2}\kappa_{\alpha\beta}}
{\partial \ary{u}^{A}\partial \ary{u}^{B}}
\Bigg]
d\bar{\Gamma}.
\label{eq:total_stiffness_matrix}
\end{align}
The first four terms in Eq.~\eqref{eq:total_stiffness_matrix} constitute the material stiffness arising from the constitutive tangent operators $\mathbb{A}$, $\mathbb{B}$ and $\mathbb{D}$, whereas the remaining terms form the geometric stiffness associated with the current stress state. These geometric contributions are essential for accurately capturing large-deformation effects, bifurcation phenomena and snap-through instabilities.

After assembly of all element contributions, the global nonlinear equilibrium equations take the form
\begin{equation}
\ary{R}(\ary{u})
=
\ary{F}^{\mathrm{int}}(\ary{u})
-
\ary{F}^{\mathrm{ext}}
=
\ary{0},
\label{eq:global_equilibrium}
\end{equation}
where $\ary{u}$ denotes the vector collecting all displacement degrees of freedom.
Linearisation of Eq.~\eqref{eq:global_equilibrium} yields the global tangent system
\begin{equation}
\ary{K}\,\varDelta\ary{u}
=
-\ary{R},
\label{eq:newton_system}
\end{equation}
with
\begin{equation}
\ary{K}
=
\ary{K}_{\mathrm{mat}}
+
\ary{K}_{\mathrm{geo}},
\end{equation}
where $\ary{K}_{\mathrm{mat}}$ and $\ary{K}_{\mathrm{geo}}$ denote the assembled material and geometric stiffness matrices, respectively.

The resulting nonlinear system is solved iteratively using a Newton--Raphson procedure. For problems involving limit points and snap-through instabilities, the Newton iterations are combined with the arc-length continuation strategy described in the following section.

\section{Bifurcation Tracking and Post-Buckling Analysis}
\label{sec:bifurcation}
Algorithm~\ref{alg:bifurcation} outlines a staged arc-length procedure combined with eigenmode perturbation for bifurcation tracking. A detailed description is given in the subsequent sections.

\begin{algorithm}[t]
\caption{Staged arc-length porcedure with eigenmode perturbation for bifurcation tracking}
\label{alg:bifurcation}
\begin{algorithmic}[0]
\Require Initialise mesh, material parameters, target loads $p_{\text{target}}$, $\Delta\Phi_{\text{target}}$, arc-length $\tilde\varDelta s$, tolerance $\varepsilon$
\Ensure Equilibrium paths $\{(\mathbf{u}_n,\mathcal{L}_n)\}_{n=1}^{N}$

\State \textbf{Stage 1: Pressure loading without electric field}
\State Set $\mathbf{u}_0 \leftarrow \mathbf{0}$, $\mathcal{L}_0 \leftarrow 0$, $p_{\text{mid}}$
\While{$p < p_{\text{mid}}$}
    \State Solve $\mathbf{R}(\mathbf{u},\mathcal{L})=\mathbf{0}$ with arc-length constraint
    \State Update $\mathbf{u}_n$, $\mathcal{L}_n$
\EndWhile
\State Store converged state $(\mathbf{u}_{\text{mid}},p_{\text{mid}})$ as initial guess for Stage~2

\State \textbf{Stage 2: Voltage loading at fixed pressure}
\State Restore $(\mathbf{u}_0,\mathcal{L}_0) \leftarrow (\mathbf{u}_{\text{mid}},p_{\text{mid}})$
\While{$\Delta\Phi < \Delta\Phi_{\text{target}}$}
    \State Solve $\mathbf{R}(\mathbf{u},\mathcal{L})=\mathbf{0}$ with arc-length constraint
    \State Update $\mathbf{u}_n$, $\mathcal{L}_n$
\EndWhile
\State Store converged state $(\mathbf{u}_{\text{elec}},\Delta\Phi_{\text{target}})$

\State \textbf{Stage 3: Further pressure loading at fixed voltage}
\State Restore $(\mathbf{u}_0,\mathcal{L}_0) \leftarrow (\mathbf{u}_{\text{elec}},\Delta\Phi_{\text{target}})$
\State Apply eigenmode perturbation to break symmetry: 
\State \quad compute eigenvector $\boldsymbol{\psi}_1$ of $\mathbf{K}(\mathbf{u}_0)$ when zero eigenvalue detected, set $\mathbf{u}_0 \leftarrow \mathbf{u}_0 + \epsilon\,\boldsymbol{\psi}_1$ with $\epsilon \ll 1$
\While{$p < p_{\text{target}}$}
    \State Solve $\mathbf{R}(\mathbf{u},\mathcal{L})=\mathbf{0}$ with arc-length constraint
    \If{new bifurcation detected}
        \State Re-apply eigenmode perturbation if needed
    \EndIf
    \State Update $\mathbf{u}_n$, $\mathcal{L}_n$
\EndWhile
\State \Return All converged equilibrium points $(\mathbf{u}_n,\mathcal{L}_n)$
\end{algorithmic}
\end{algorithm}
\subsection{Staged Arc-Length Procedure}

To capture symmetry-breaking bifurcations and post-buckling equilibrium paths, the nonlinear equilibrium equations are solved using a Newton--Raphson scheme combined with a Crisfield spherical arc-length method. 
The arc-length strategy enables the solution path to be traced through limit points and snap-through instabilities that cannot be followed by conventional load-control procedures.

The external loading is parameterized by a scalar load multiplier
\begin{equation}
\ary{F}^{\mathrm{ext}}(\mathcal{L})
=
\mathcal{L}\,\overline{\ary{F}},
\end{equation}
where $\mathcal{L}$ denotes the active loading parameter. Depending on the loading stage, $\mathcal{L}$ represents either the pressure scaling factor
\(
p=\mathcal{L}p_0
\)
or the voltage scaling factor
\(
\Delta\Phi=\mathcal{L}\Delta\Phi_0.
\)

The nonlinear equilibrium condition is expressed as
\begin{equation}
\ary{R}(\ary{u},\mathcal{L})
=
\ary{F}^{\mathrm{int}}(\ary{u})
-
\ary{F}^{\mathrm{ext}}(\mathcal{L})
=
\ary{0}.
\end{equation}

Since both the displacement field $\ary{u}$ and the load multiplier $\mathcal{L}$ are unknown, an additional arc-length constraint is introduced,
\begin{equation}
g(\tilde\varDelta\ary{u},\tilde\varDelta\mathcal{L})
=
\tilde\varDelta\ary{u}^{\mathsf T}\tilde\varDelta\ary{u}
+
\psi^2(\tilde\varDelta\mathcal{L})^2
-
\tilde\varDelta s^2
=
0,
\end{equation}
where $\tilde\varDelta s$ is the prescribed arc length and $\psi$ is a scaling parameter controlling the relative contribution of displacement and load increments.

Here, $\tilde\delta(\cdot)$ denotes the Newton correction within the current iteration, whereas $\tilde\varDelta(\cdot)$ denotes the cumulative increment over the current arc-length step.

At Newton iteration $k$, the equilibrium equations and the arc-length constraint are linearized simultaneously, resulting in the extended system
\begin{equation}
\begin{bmatrix}
\ary{K}^{(k)}
&
-\overline{\ary{F}}
\\
2\tilde\varDelta\ary{u}^{(k)\mathsf T}
&
2\psi^2\tilde\varDelta\mathcal{L}^{(k)}
\end{bmatrix}
\begin{bmatrix}
\tilde\delta\ary{u}
\\
\tilde\delta\mathcal{L}
\end{bmatrix}
=
-
\begin{bmatrix}
\ary{R}^{(k)}
\\
g^{(k)}
\end{bmatrix},
\end{equation}
where $\ary{K}^{(k)}$ is the consistent tangent stiffness matrix and
\begin{equation}
g^{(k)}
=
\tilde\varDelta\ary{u}^{(k)\mathsf T}
\tilde\varDelta\ary{u}^{(k)}
+
\psi^2
\left(
\tilde\varDelta\mathcal{L}^{(k)}
\right)^2
-
\tilde\varDelta s^2
\end{equation}
is the residual of the arc-length constraint.

After solving for the Newton corrections, the cumulative increments are updated as
\begin{equation}
\tilde\varDelta\ary{u}^{(k+1)}
=
\tilde\varDelta\ary{u}^{(k)}
+
\tilde\delta\ary{u},
\qquad
\tilde\varDelta\mathcal{L}^{(k+1)}
=
\tilde\varDelta\mathcal{L}^{(k)}
+
\tilde\delta\mathcal{L}.
\end{equation}

The total solution is then updated according to
\begin{equation}
\ary{u}^{(k+1)}
=
\ary{u}_{n}
+
\tilde\varDelta\ary{u}^{(k+1)},
\qquad
\mathcal{L}^{(k+1)}
=
\mathcal{L}_{n}
+
\tilde\varDelta\mathcal{L}^{(k+1)},
\end{equation}
where $(\ary{u}_{n},\mathcal{L}_{n})$ denotes the converged solution at the previous arc-length step.

The Newton iterations continue until both the residual norm and the displacement increment satisfy the prescribed convergence tolerances.

To reproduce the experimental loading protocol and improve numerical robustness, the equilibrium path is traced through three consecutive continuation stages:
\begin{enumerate}
\item Pressure loading from the undeformed configuration to an intermediate pressure level $p_{\mathrm{mid}}$ with $\Delta\Phi=0$;

\item Voltage loading from $0$ to $\Delta\Phi_{\mathrm{target}}$ while maintaining the pressure fixed at $p_{\mathrm{mid}}$;

\item Further pressure loading from $p_{\mathrm{mid}}$ to the target pressure while keeping
$\Delta\Phi=\Delta\Phi_{\mathrm{target}}$ constant.
\end{enumerate}

The converged solution of each stage serves as the initial state for the subsequent stage. This staged arc-length procedure allows the complete electro-mechanical equilibrium manifold to be followed, including symmetry-breaking bifurcations, limit points, and snap-through transitions.

\subsection{Mechanism of Inducing Bifurcation}
The electromechanical response of a perfectly symmetric shell under internal pressure and voltage loading is inherently susceptible to symmetry-breaking bifurcations. Due to the electro-softening effect, the membrane's stiffness degrades as the electric potential increases. Upon reaching a critical load, the axisymmetric equilibrium path loses stability, and the structure seeks a lower-energy configuration by transitioning to a non-axisymmetric deformed state.
Standard incremental solvers fail at the bifurcation point because the perfectly symmetric discretization lacks the necessary perturbation to trigger the descent onto the asymmetric branch. Several techniques exist to overcome this, including eigenmode perturbation, geometric imperfections, and material imperfections. In this work, we employ a specific numerical strategy to navigate this singularity and capture the post-buckling response.

{
\textbf{Remarks:} The perturbation is scaled to a small fraction of the shell's characteristic dimension. However, determining the optimal scale factor is non-trivial: if too small, the solver may remain trapped on the principal branch; if too large, it may overshoot the true post-buckling path. In practice, the scale requires careful calibration based on mesh resolution and problem-specific sensitivity to ensure robust convergence onto the asymmetric branch without artificially altering the physical response.
}
\section{Numerical Examples}
\label{sec:numerical_examples}
This section presents three numerical examples designed to systematically validate and demonstrate the capabilities of the proposed electroelastic thin shell formulation. 
The examples are sequenced to progress from fundamental verification to the exploration of complex, coupled nonlinear phenomena. 
First, the inflation of a spherical membrane serves to benchmark the model’s accuracy against a known analytical solution. 
Second, the electromechanical response of a prestretched circular plate under combined pressure and voltage loading is investigated to elucidate the competitive interplay between mechanical prestress and electric-field-induced softening. 
Finally, the primary focus is on the nonlinear buckling and symmetry-breaking behaviour of toroidal membranes, showcasing the formulation’s ability to trace complex post-bifurcation equilibrium paths under strong electromechanical coupling. 
All simulations are performed using the numerical implementation based on subdivision surfaces described in Section~\ref{sec:implementation}.
\subsection{Validation: Electroelastic Spherical Shells}
\begin{figure}[h]
\centering
\includegraphics[width=\linewidth]{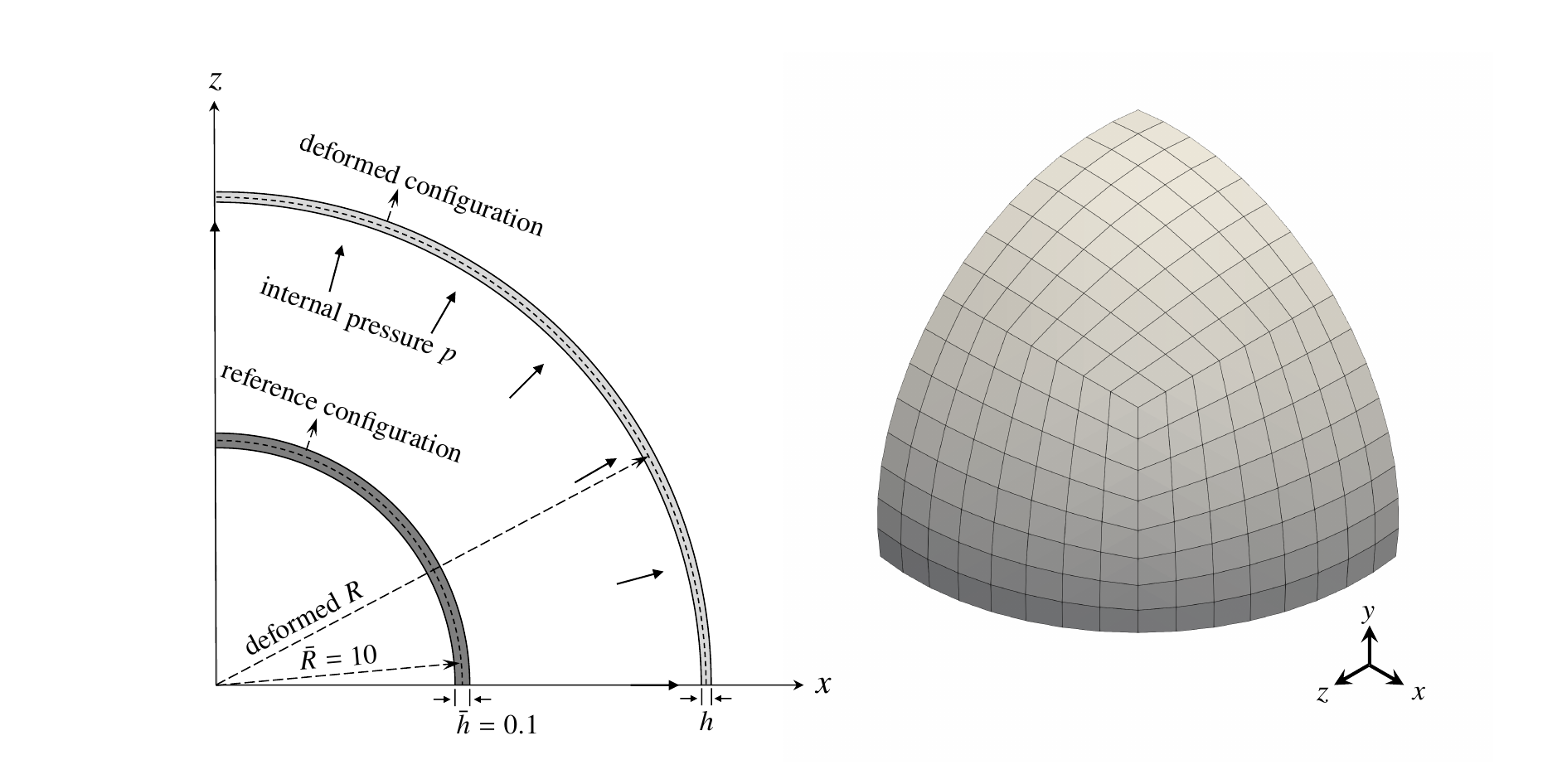}
\caption{Geometric definition of the electroelastic spherical shell. Left: reference and deformed configurations (slice on the $xz$-plane). Right: control mesh representing only one-eighth of the hemisphere in reference configuration, utilising symmetry boundary conditions for the analysis.}
\label{fig:sphere_geometry}
\end{figure}
\begin{figure}[h]
\centering
\includegraphics[width=0.8\linewidth]{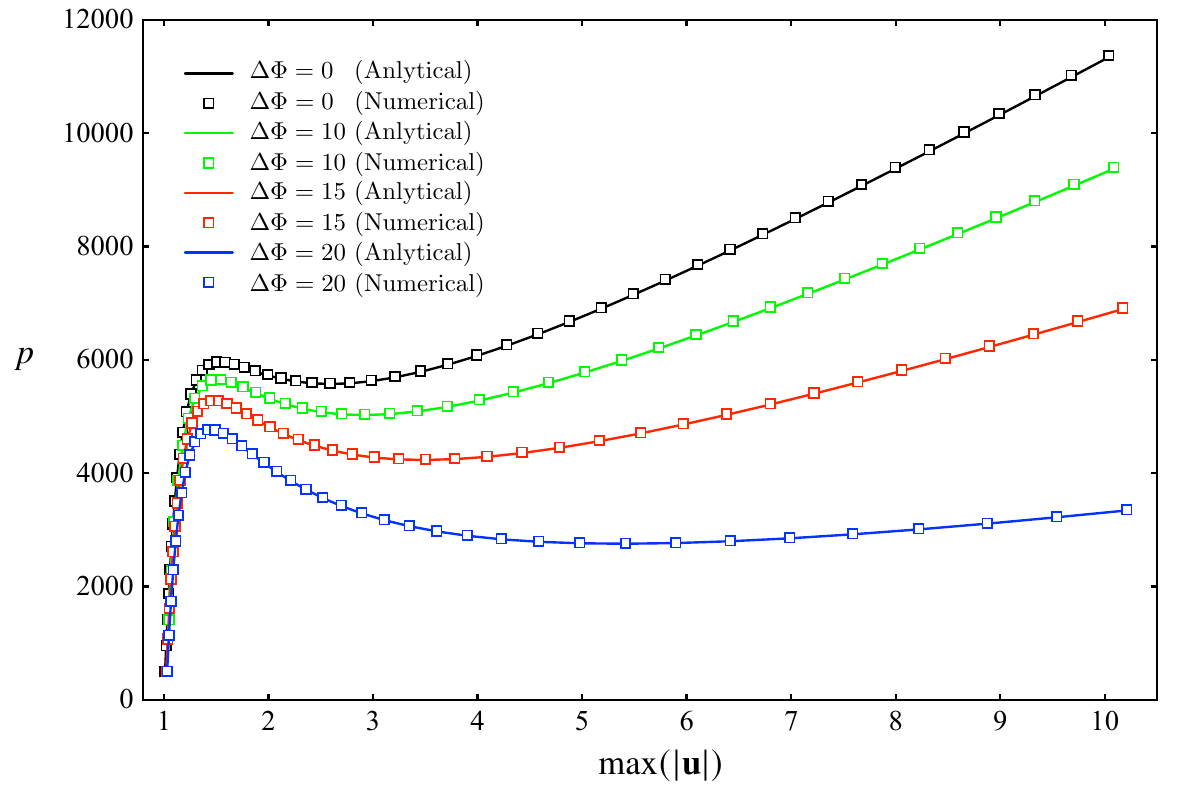}
\caption{Variation of inflation pressure with the maximum displacement for different electrical potentials. }
\label{fig:elec-sphere}
\end{figure}
To validate the proposed method, the inflation behaviour of an electroelastic spherical shell is examined against available analytical solutions. 
The pure-mechanical case serves as a widely adopted benchmark, as demonstrated in our previous work~\cite{liu2024computational} and by others~\cite{chen2014explicit, Cirak:2001aa, kiendl2015isogeometric}. 
The analytical electroelastic solution employed here is derived by extending this benchmark to incorporate electromechanical coupling effects.
The geometric setup and the corresponding control mesh are presented inFig.~\ref{fig:sphere_geometry}.
One considers an electroelastic spherical shell of radius $\bar{R} = 10$ and thickness $\bar{h} = 0.1$, subjected to both internal pressure and an electrical potential. An inflating pressure $p$ is applied to the inner surface of the shell, while different electric potentials $\Delta\Phi \in \{0, 10, 15, 20\}$ are individually imposed across its thickness. The material parameters are taken as $c_1 = 0.4375\mu$ and $c_2 = 0.0625\mu$, while $\mu = 4.225\times10^5$. 
For the spherical shell, the stretching of the mid-surface is the same in all directions, thus $\lambda_1 = \lambda_2 = \lambda$.
The analytical solution for the inner pressure ${p}$ is given by:
\begin{equation}
p = \frac{4\bar{h}}{\bar{R}} \left[ c_1 [\lambda^{-1} - \lambda^{-7}] - c_2 [\lambda^{-5} - \lambda] - \frac{\varepsilon [\Delta\Phi]^2 \lambda}{2 \bar{h}^2} \right].
\end{equation}
Fig.~\ref{fig:elec-sphere} compares the numerical results obtained from the proposed method with the analytical solution. The plot shows the relationship between the internal pressure $p$ and the maximum radial displacement $u$ under different applied electric potentials. The solid lines represent the analytical prediction derived from the variational principle, while the markers denote the numerical results from the present finite element implementation.
As shown in the figure, the numerical results are in excellent agreement with the analytical curves across all voltage levels. The electrostatic softening effect, manifested as a leftward and downward shift of the pressure–displacement curve with increasing $\Delta\Phi$, is accurately captured by the numerical model. This validation confirms the correctness and robustness of the proposed electroelastic thin shell formulation for large-deformation problems involving strong electromechanical coupling.

\subsection{Electro-Softening: Inflation of Prestreched Circular Plates}
\begin{figure}[h]
\centering
\includegraphics[width=1.05\linewidth]{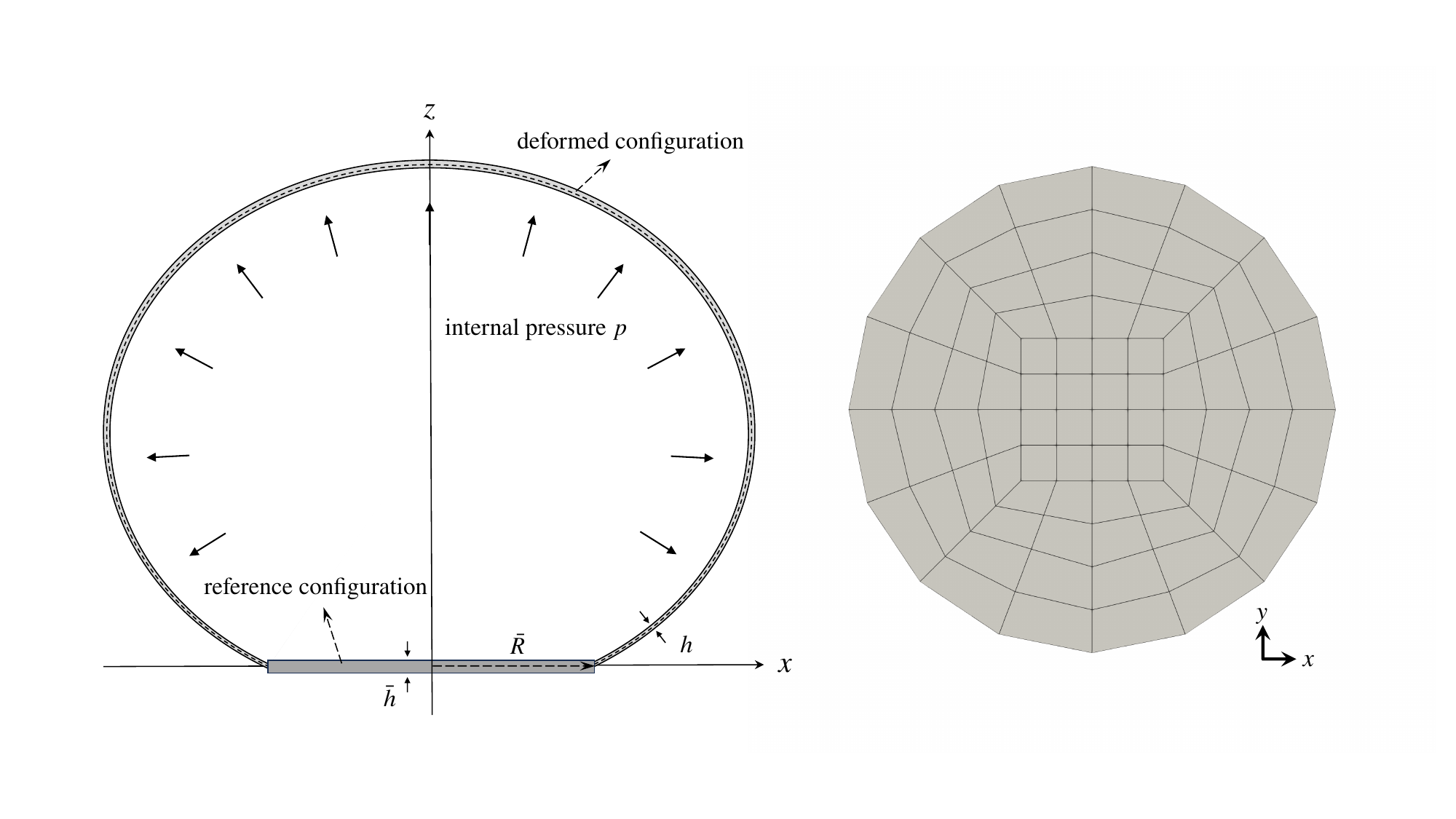}
\caption{
Geometric definition of the electroelastic circular plate. Left: reference and deformed configurations (slice on the $xz$-plane). Right: control mesh for reference configuration.
}
\label{fig:plate_geometry}
\end{figure}
This example investigates the coupled response of a prestretched dielectric elastomer membrane subjected to simultaneous inflation pressure and electric potential (shown in Fig.~\ref{fig:plate_geometry}). 
The primary objective is to analyse the competing mechanisms between mechanical prestretch, which increases the effective stiffness, and the applied electric field, which induces a pronounced electro-softening effect. 
A circular dielectric elastomer membrane with an initial undeformed radius $\bar{R} = 7.5 \times 10^{-2}$ and thickness $\bar{h} = 5\times10^{-4}$ is considered. 
The material parameters for the Mooney–Rivlin model are $c_1 = 8\times10^4$, $c_2 = 2\times10^4$.

\begin{figure}[h]
\centering
\begin{subfigure}[b]{0.48\linewidth}
  \centering
  \includegraphics[width=\linewidth, height=0.25\textheight, keepaspectratio]{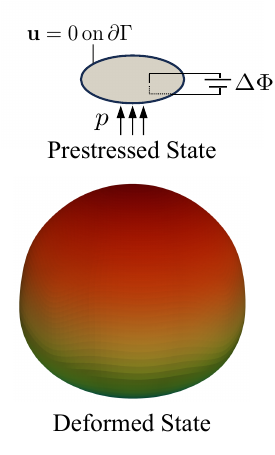}
  \caption{Deformed and undeformed states}
  \label{fig:Deformation of the plate}
\end{subfigure}
\hfill
\begin{subfigure}[b]{0.48\linewidth}
  \centering
  \includegraphics[width=\linewidth]{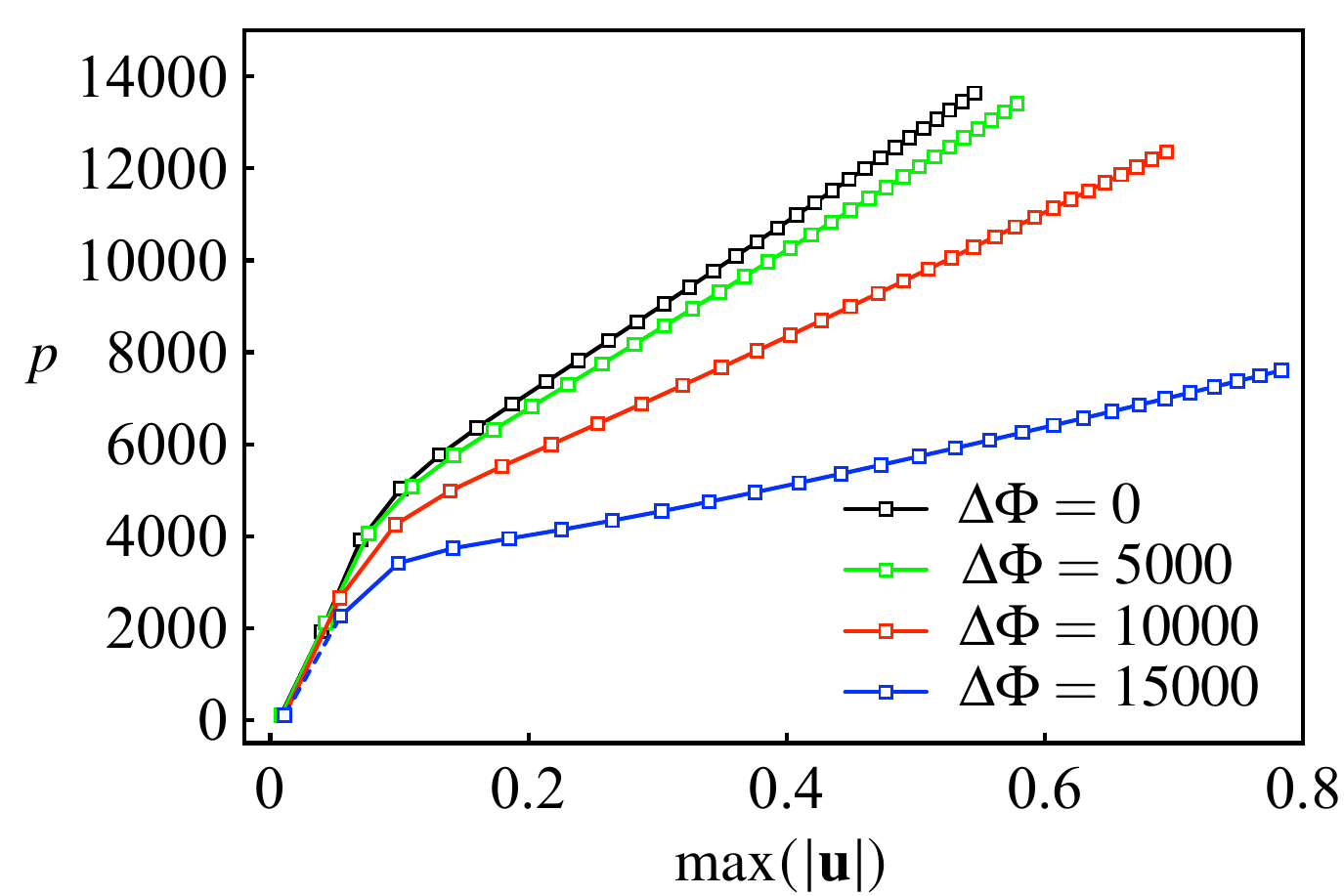}
  \caption{ $\lambda_p = 1.1$}
  \label{fig:lambda_p=1.1}
\end{subfigure}
\vspace{0.1cm} 
\begin{subfigure}[b]{0.48\linewidth}
  \centering
  \includegraphics[width=\linewidth]{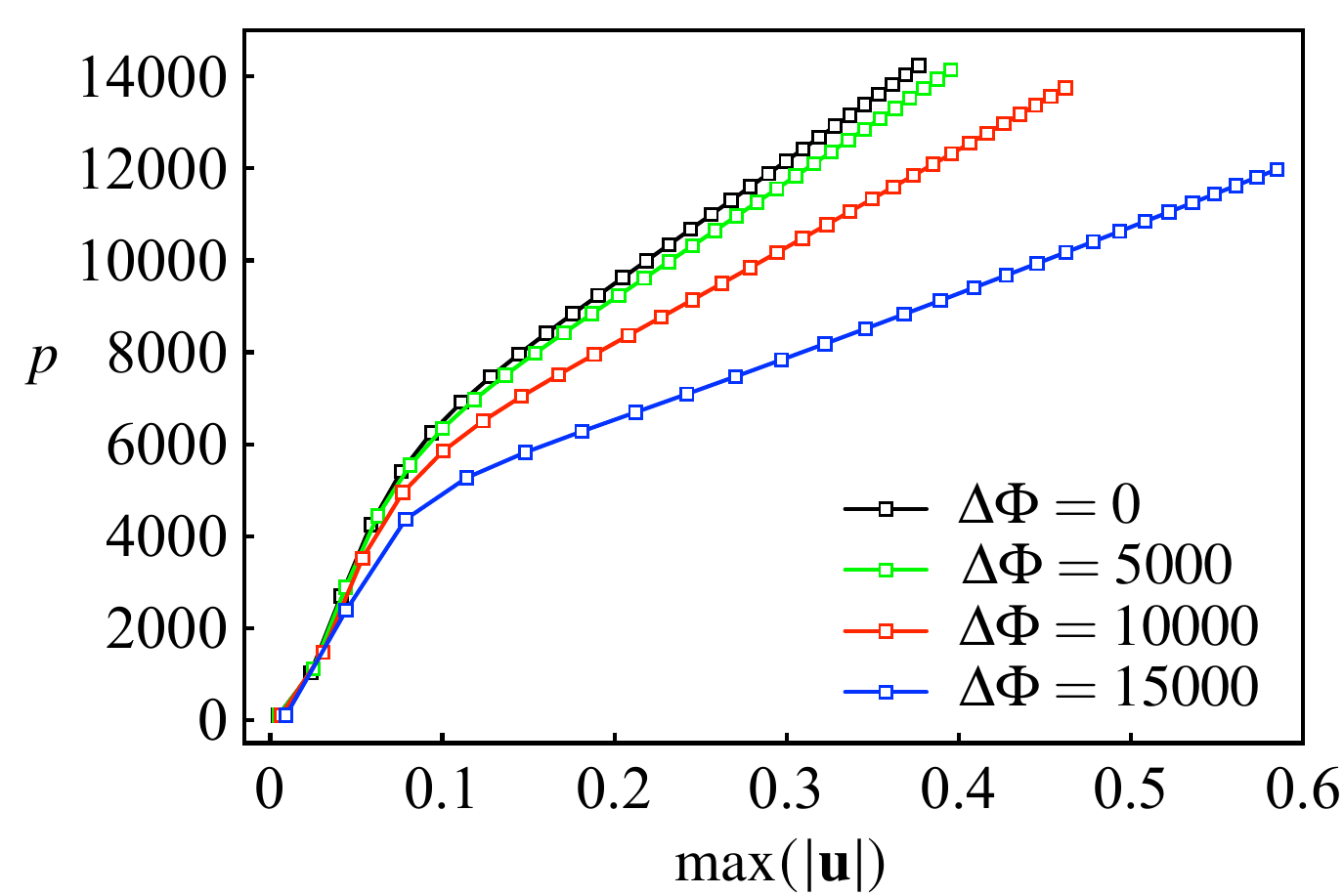}
  \caption{ $\lambda_p = 1.2$}
  \label{fig:lambda_p=1.2}
\end{subfigure}
\hfill
\begin{subfigure}[b]{0.48\linewidth}
  \centering
  \includegraphics[width=\linewidth]{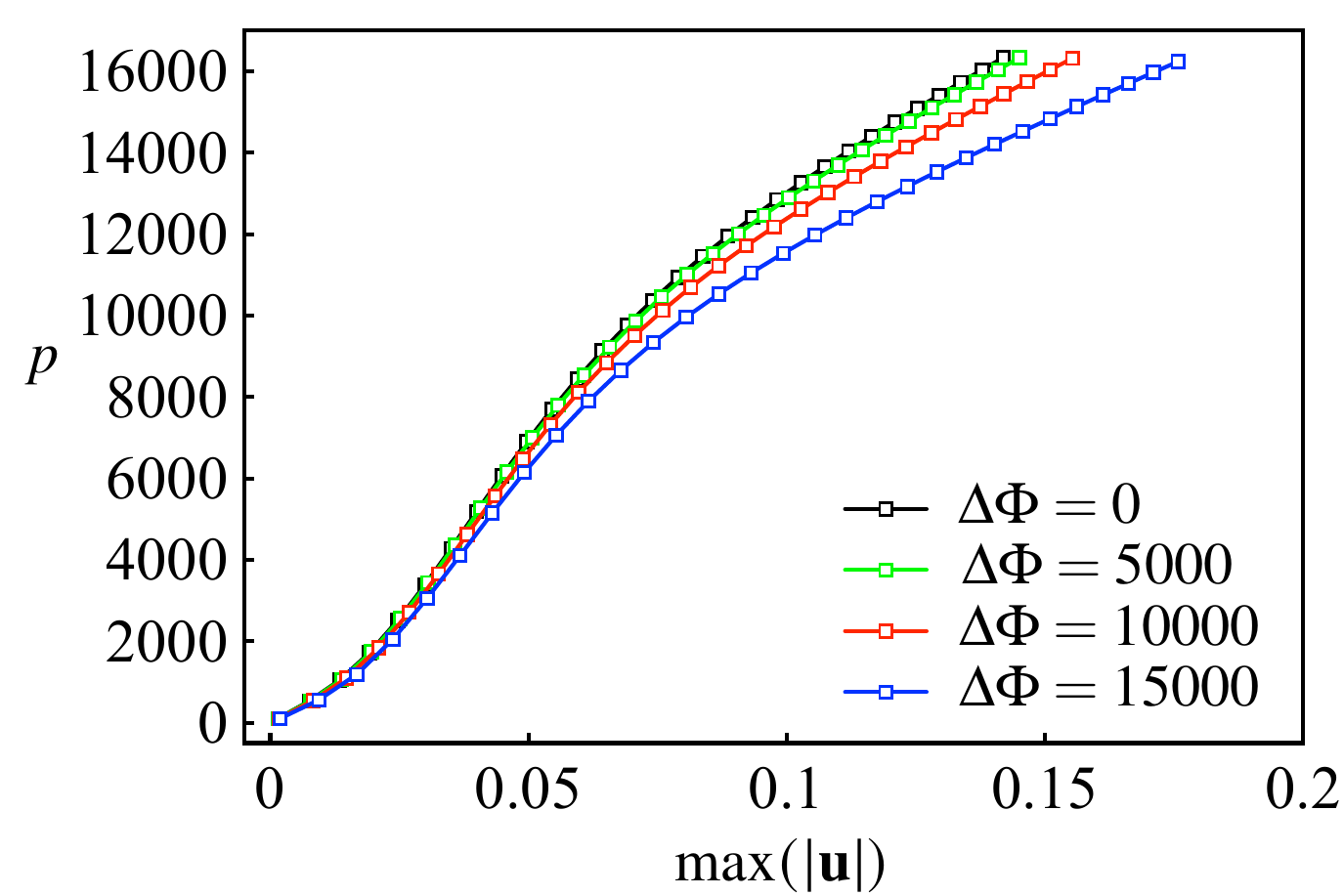}
  \caption{ $\lambda_p = 1.5$}
  \label{fig:lambda_p=1.5}
\end{subfigure}
\caption{Variation curves of internal pressure $p$ versus maximum displacement for a prestreched dielectric elastomer circular membrane under electrical loading, with electric potential $\Delta\Phi$ taking values of \{{0}, {5000}, {10000}, {15000}\}.}
\label{fig:plate}
\end{figure}

Fig.~\ref{fig:plate} illustrates the load–displacement curves and corresponding deformation profiles of the electro-elastic circular membrane under various levels of prestrech.
The electric potential \(\Delta\Phi\) increases from {0} to {15000} with a step size of {5000}, while the stretch ratio \(\lambda_p\) induced by the prestress varies within the set \(\{1.1, 1.2, 1.5\}\).
The pressure–displacement curve flattens as the electric potential rises, indicating a loss of structural stiffness. Accordingly, the maximum displacement at a constant internal pressure increases substantially.
Introducing a prestrech elevates the membrane's load-bearing capacity by enhancing its effective stiffness. 
This increased structural rigidity results in a substantial suppression of the maximum displacement.
As illustrated in Fig.~\ref{fig:lambda_p=1.1}, upon applying an electric potential of {15000}, the maximum displacement under a given internal pressure exhibits a marked surge compared to the zero-potential scenario. 
This behaviour stems from the degradation of material stiffness induced by the electro-softening effect.
This is attributed to the reduction in material stiffness caused by the electro-softening effect.
In contrast, in Fig.~\ref{fig:lambda_p=1.5}, under the same applied potential of {15000}, the increase in the maximum displacement for the same internal pressure is relatively small. 
This indicates that, for the present numerical example, an increase in prestress suppresses the electro-softening effect. 
\subsection{Symmetric Breaking of Electroelastic Toroidal Membranes}
The symmetry-breaking behaviour of electroelastic toroidal membranes was investigated in our previous work~\cite{liu2020coupled} using a semi-analytical framework. 
Subsequently, a simplified discrete model for axisymmetric dielectric elastomer membranes was developed~\cite{liu2024simplified} and successfully applied to the analysis of toroidal structures.
Although the axisymmetric formulation greatly reduces the computational cost by exploiting geometric symmetry, it inherently excludes non-axisymmetric deformation modes.
Consequently, neither symmetry-breaking bifurcations nor the associated post-buckling branches can be resolved.
The present three-dimensional electroelastic shell formulation overcomes this limitation, providing a general framework for analysing symmetry-breaking instabilities, branch switching, and post-buckling behaviour in toroidal dielectric elastomer structures.
\begin{figure}[h]
\centering
\includegraphics[width=1.05\linewidth]{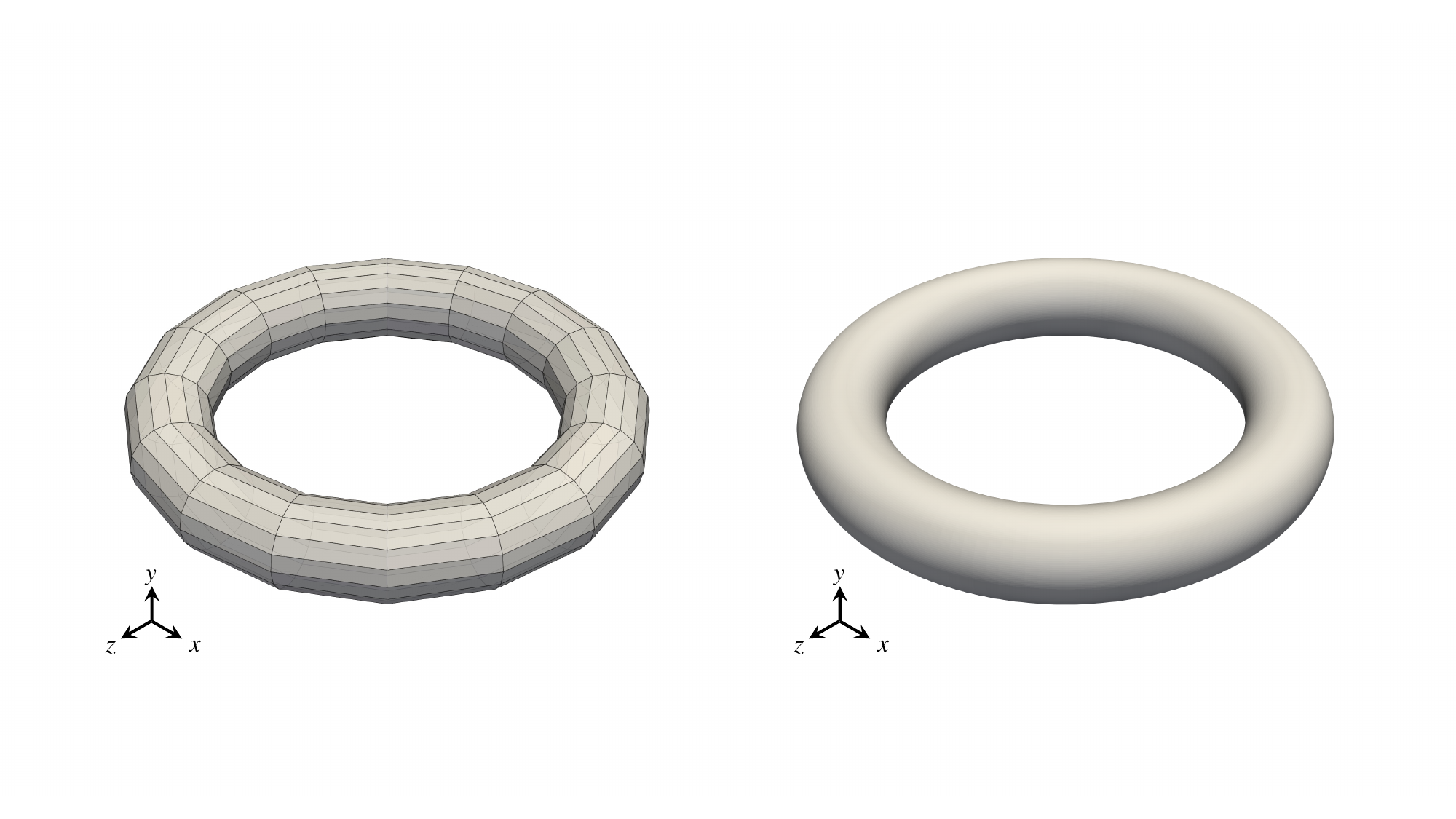}
\caption{The left figure shows that in the numerical analysis, the toroidal shell is discretised into $256$ elements with a total of $256$ nodes; the right figure illustrates the smooth limiting surface of the toroidal shell.
}
\label{fig:torus_geometry}
\end{figure}
\begin{figure}[h]
\centering
\begin{subfigure}[b]{0.48\linewidth}
  \includegraphics[width=\linewidth]{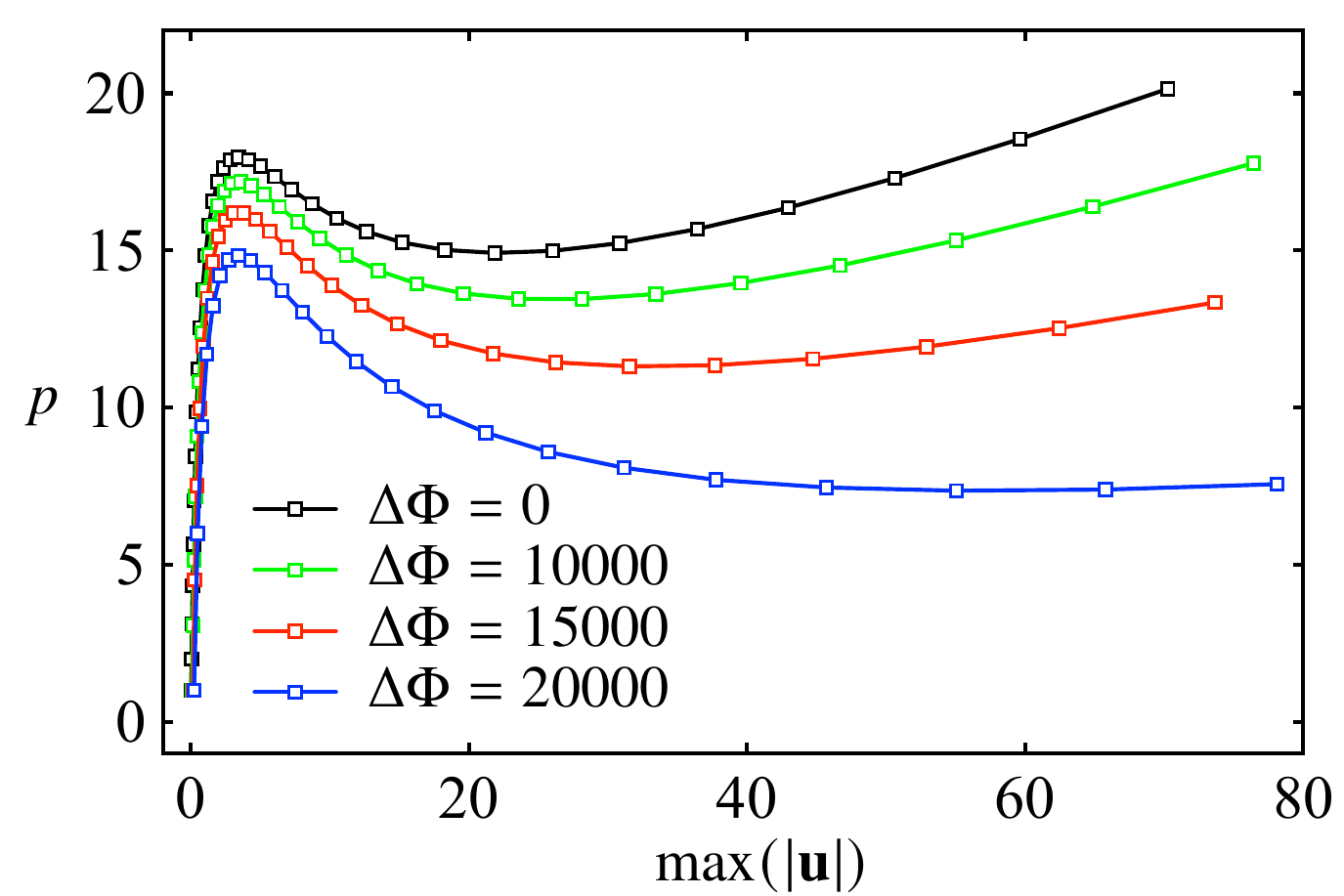}
  \caption{}
  \label{fig:elec_torus_0.2}
\end{subfigure}
\hfill
\begin{subfigure}[b]{0.48\linewidth}
  \includegraphics[width=\linewidth]{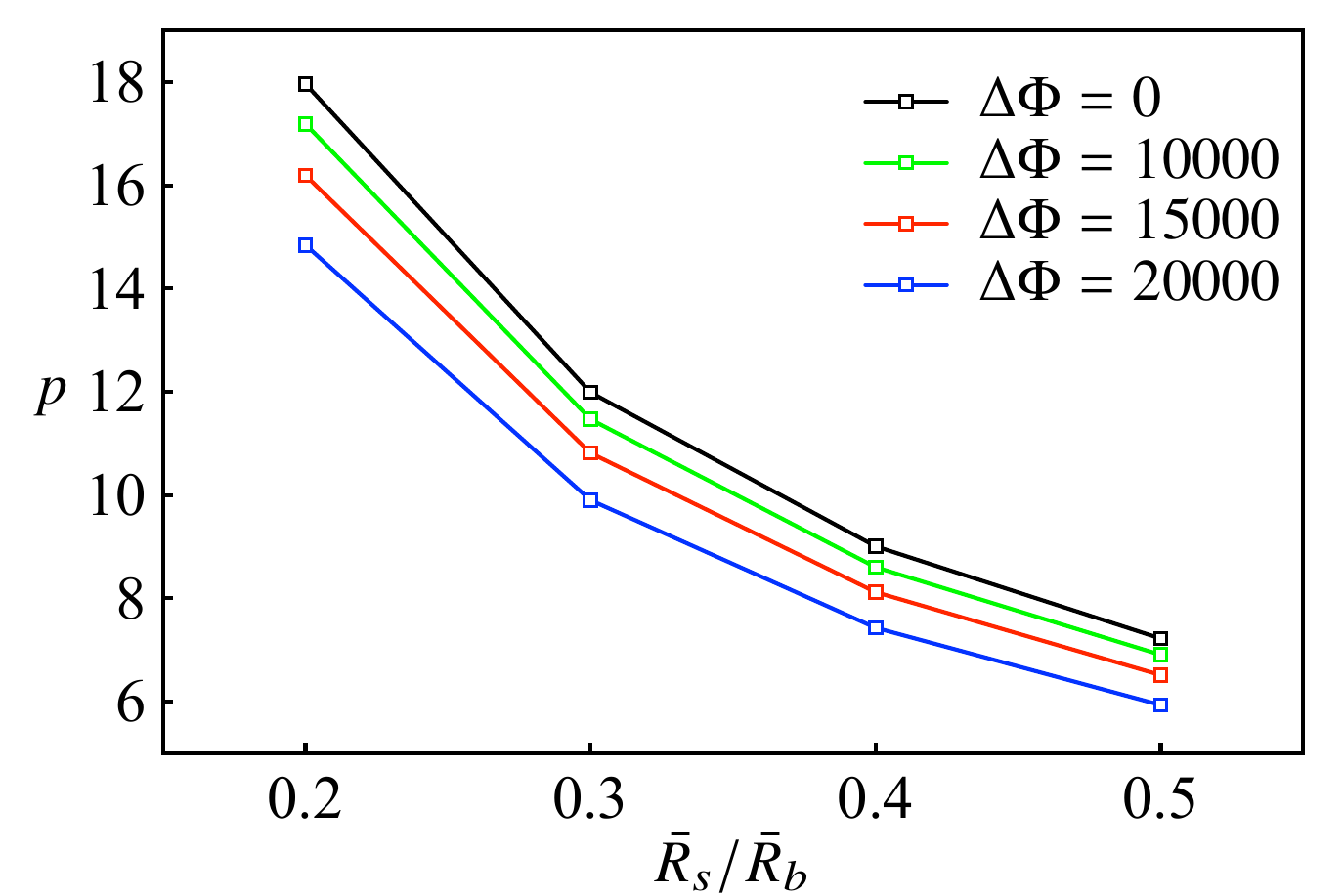}
  \caption{}
  \label{fig:elec_torus_aspect_ratio}
\end{subfigure}
\caption{ (a) Variation curves of the pressure with respect to the maximum displacement for a toroidal membrane with torus radius $\bar{R}_b = 10$, cross-sectional radius $\bar{R}_s = 2$, and thickness $\bar{h} = 0.01 $ under electrical loading, where the electric potential $\Delta\Phi$ varies within the set $\{0, 10000, 15000, 20000\}$; (b) Relationship between the limit point pressure $p_s$ and the aspect ratio $\bar{R}_s/\bar{R}_b$ under electrical loading, with the electric potential $\Delta\Phi$ varying within the set $\{0, 10000, 15000, 20000\}$. }
\label{fig:two_images}
\end{figure}
\subsubsection{Symmetric Inflation and Limit Point Instability}
We first investigate the deformation behaviour and principal equilibrium pathways under varying electric fields. 
A numerical example is first established to simulate the inflation process, employing material parameters $c_1 = 0.4375\mu$ and $c_2 = 0.0625\mu$ for a torus with major radius $\bar{R}_b = 10$, minor radius $\bar{R}_s = 2$, and thickness $\bar{h} = 0.01$.
The control mesh and limiting surface are shown in Fig.~\ref{fig:torus_geometry}.

Figure~\ref{fig:elec_torus_0.2} illustrates the variation of the pressure $p$ with respect to the maximum displacement $u_m$ for different electrical loads. 
The results indicate that the application of an electric load effectively reduces the load-bearing capacity of the structure.

From a stability perspective, the limit point corresponds to the loss of tangent stiffness and the onset of snap-through instability. 
Under pure mechanical loading, the limit pressure $p_s$ marks the transition from a stable equilibrium branch to an unstable post-buckling path. 
When an electric field is applied, the Maxwell stress introduces an additional tensile component in the membrane, which reduces the effective structural stiffness and shifts the equilibrium path towards lower pressures. 
Consequently, the limit point moves to smaller displacements and lower critical loads, indicating that the electric field acts as a destabilising agent.

Furthermore, Figure~\ref{fig:elec_torus_aspect_ratio} depicts the variation of the limit-point pressure $p_s$ with the aspect ratio $\bar{R}_s/\bar{R}_b$. 
It is observed that $p_s$ decreases monotonically with increasing aspect ratio and that higher electrical loads further suppress the critical buckling pressure. 
These numerical observations are in excellent agreement with existing analytical benchmarks~\cite{liu2024simplified}.
\subsubsection{Electromechanical Softening and the Role of Material Nonlinearity}
\begin{figure}[h]
\centering
\begin{subfigure}[b]{0.48\linewidth}
  \centering
  \includegraphics[width=\linewidth]{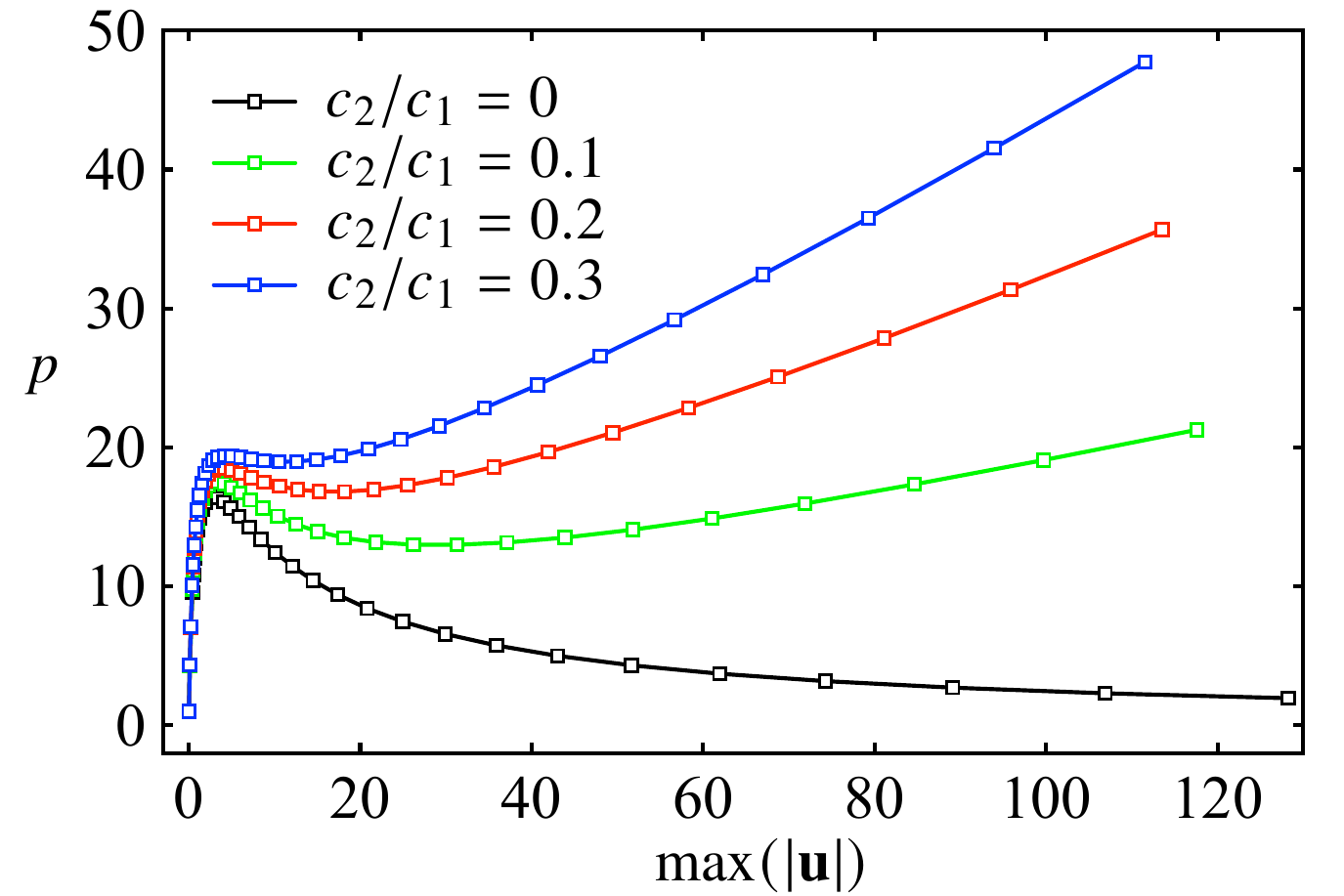}
  \caption{$\varDelta \Phi = 0$}
  \label{fig:E=0}
\end{subfigure}
\hfill
\begin{subfigure}[b]{0.48\linewidth}
  \centering
  \includegraphics[width=\linewidth]{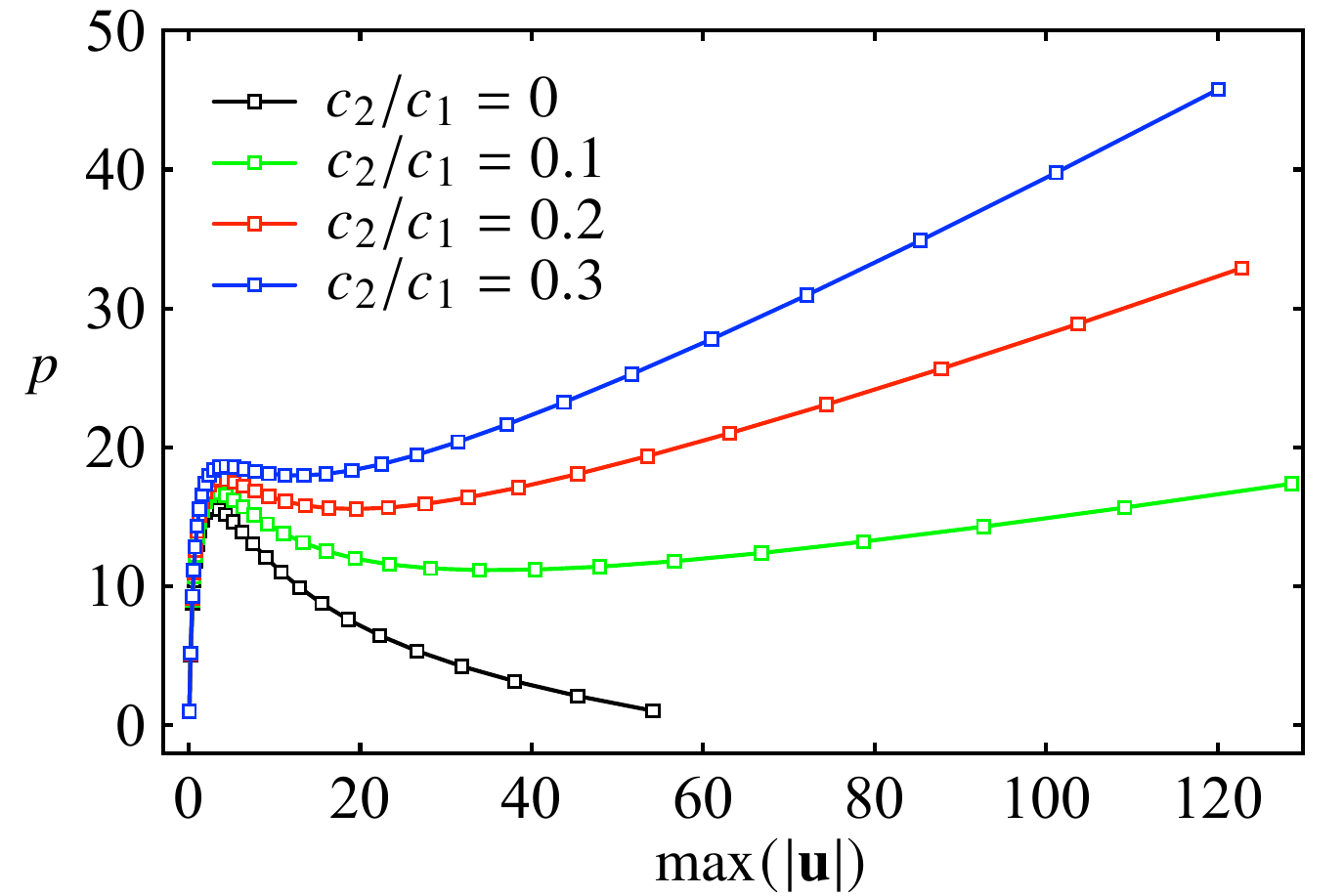}
  \caption{$\varDelta \Phi = 10000$}
  \label{fig:E=0.1538}
\end{subfigure}
\vspace{0.1cm} 
\begin{subfigure}[b]{0.48\linewidth}
  \centering
  \includegraphics[width=\linewidth]{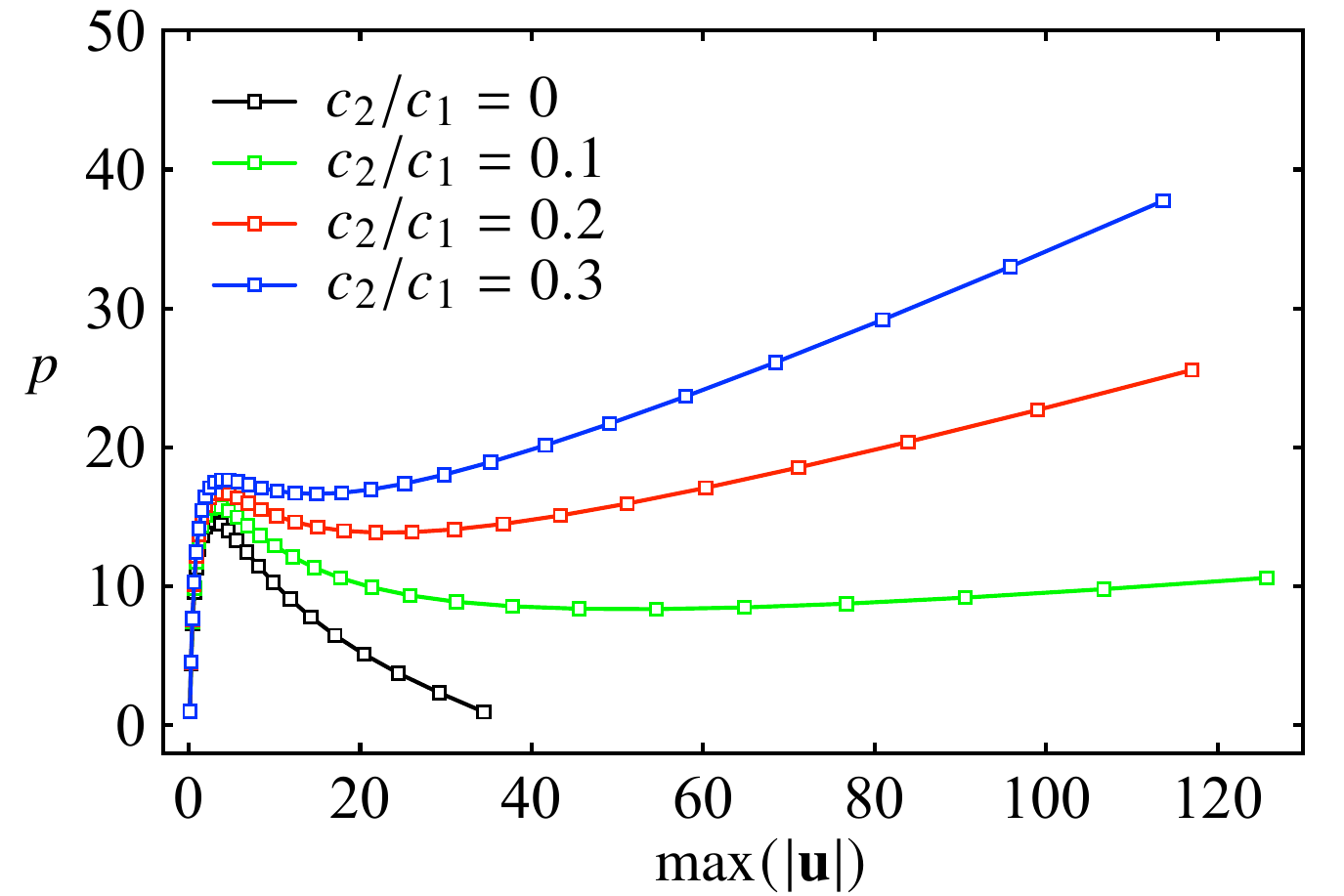}
  \caption{$\varDelta \Phi = 15000$}
  \label{fig:E=0.2308}
\end{subfigure}
\hfill
\begin{subfigure}[b]{0.48\linewidth}
  \centering
  \includegraphics[width=\linewidth]{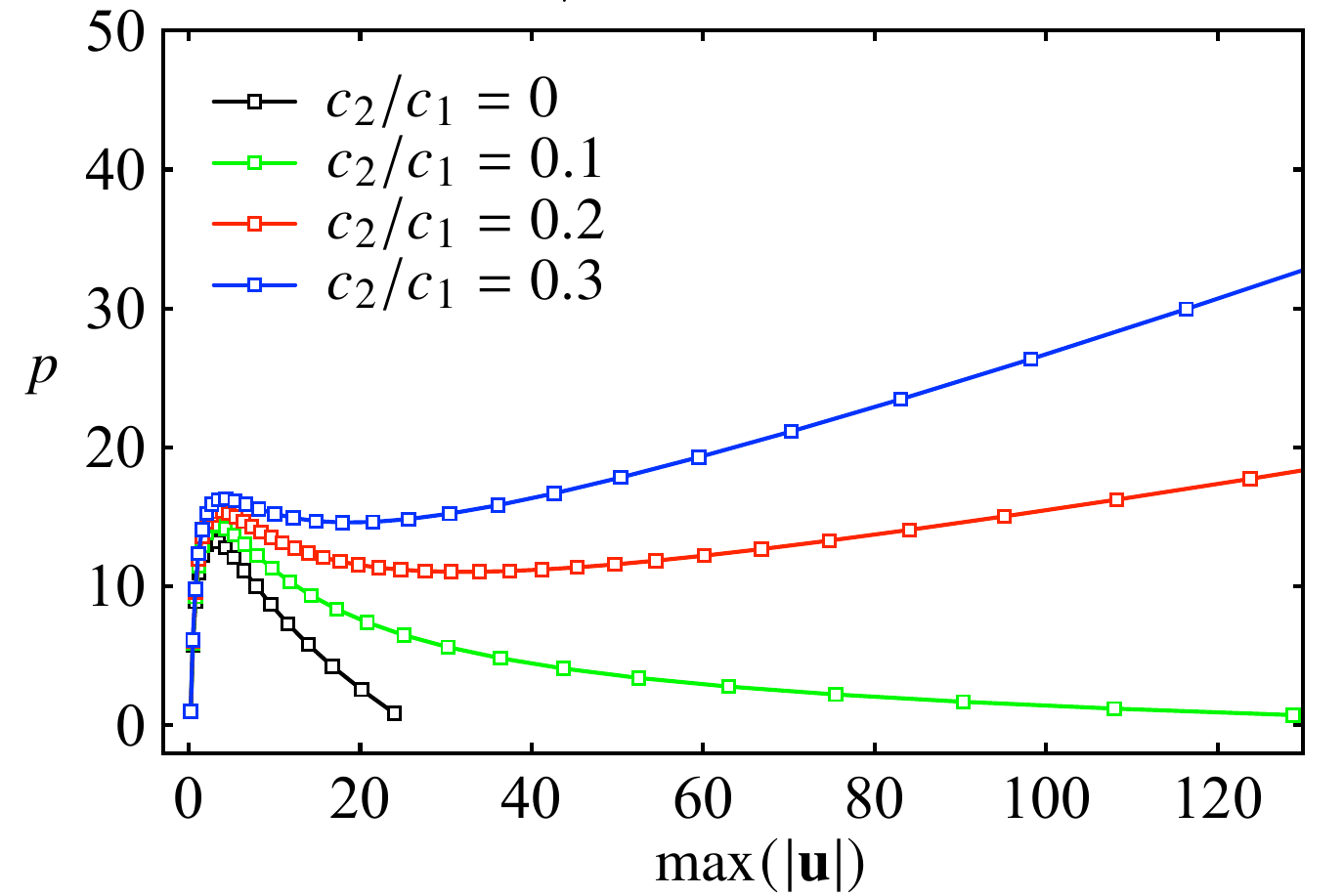}
  \caption{$\varDelta \Phi = 20000$}
  \label{fig:E=0.3077}
\end{subfigure}
\caption{Illustrates the relationship between the pressure and the displacement for an annular membrane with a toroidal radius $\bar{R} = 10$, cross-sectional radius $\bar{R}_s = 2$, and thickness $\bar{h} = 0.01$ under electrical loading. The material constant ratio $c_2/c_1$ takes values of $\{0, 0.1, 0.2, 0.3\}$. The electric potential $\Delta\Phi$ is set to $\{0, 10000, 15000, 20000\}$.}
\label{fig:four_images}
\end{figure}

The influence of material nonlinearity on the electromechanical response is investigated by varying the Mooney--Rivlin parameter ratio $c_2/c_1 \in \{0, 0.1, 0.2, 0.3\}$ under electric potentials $\Delta\Phi \in \{0, 10000, 15000, 20000\}$. 
Figure~\ref{fig:four_images} illustrates the resulting pressure--displacement curves.
In the absence of an electric field (Fig.~\ref{fig:E=0}), an increase in $c_2/c_1$ enhances the post-limit-point stiffening, indicative of a stronger strain-hardening response at large deformations. 
However, the application of an electric load exerts a destabilising effect, progressively suppressing this strain-hardening behaviour. 
For instance, at $c_2/c_1 = 0.1$, increasing $\Delta\Phi$ from $0$ to $20000$ gradually erodes the material's hardening capacity until it vanishes entirely.

\subsubsection{Electrically Induced Localisation and Post‑Buckling Path Switching}
\begin{figure}[t]
\centering
\begin{subfigure}[b]{0.8\linewidth} 
  \smash{\raisebox{0.85\height}{\makebox[0pt][l]{(a)}}}%
  \includegraphics[width=\linewidth]{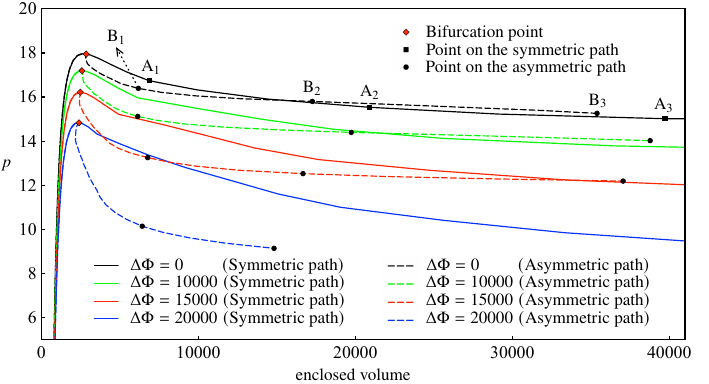}
  \label{fig:0.2_q_1}
\end{subfigure}
\hfill
\begin{subfigure}[b]{0.6\linewidth}
  \smash{\raisebox{0.5\height}{\makebox[0pt][l]{(b)}}}%
  \includegraphics[width=\linewidth]{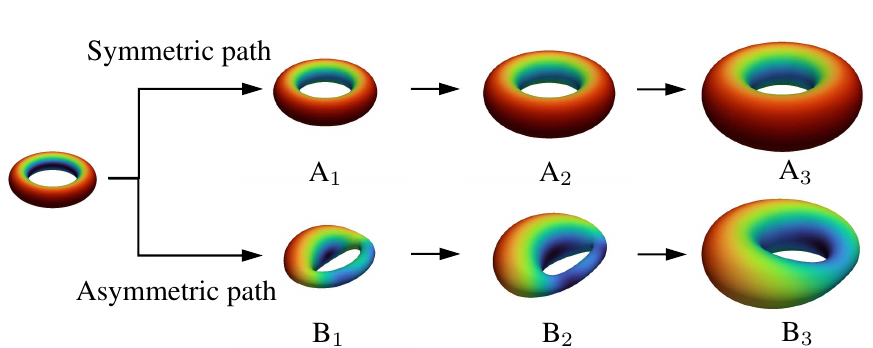}
  \label{fig:0.2_q_2}
\end{subfigure}
\begin{subfigure}[b]{0.8\linewidth}
\centering
  \smash{\raisebox{0.5\height}{\makebox[0pt][l]{(b)}}}%
  \includegraphics[width=\linewidth]{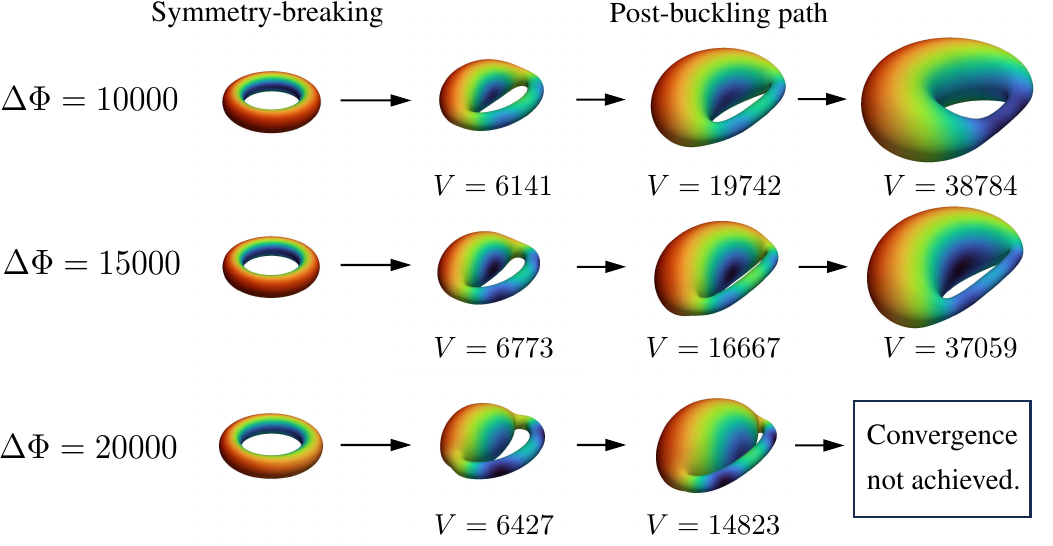}
  \label{fig:0.2_q_3}
\end{subfigure}
\caption{Pressure-volume curves and representative deformed shapes for an toroidal membrane under combined loading. The middle inset illustrates the post-buckling path switching. Bottom panels show the effect of increasing electric potential: higher voltages induce severe localised bulging and self-contact, which suppresses the symmetric equilibrium branch. The torus has $\bar{R}_b = 10$ and $\bar{R}_s = 2$.}
\label{fig:0.2_q}
\end{figure}
Tracing the post-buckling equilibrium path of a geometrically perfect toroidal membrane beyond the bifurcation point presents a significant computational challenge due to the loss of uniqueness at the critical state. 
To circumvent this numerical instability, a symmetry-breaking perturbation, either in the form of an eigenmode or a geometric imperfection, is strategically imposed.
This artificial imperfection serves as a trigger, stabilising the nonlinear solver and enabling controlled access to the post-bifurcation branches.

Figures~\ref{fig:0.2_q} and \ref{fig:0.4_q} depict the pressure--volume response alongside the corresponding deformation modes. 
The bifurcation points are observed to emerge in proximity to the limit point, with the post-buckling trajectories intersecting the principal equilibrium path. 
For the slender torus ($\bar{R}_s/\bar{R}_b = 0.2$) under high electric potential ($\Delta\Phi = 20000$), intersection is suppressed due to severe localisation: excessive deformation on one side induces self-contact, preventing further volume expansion along the symmetric branch. 
Conversely, for the thicker torus ($\bar{R}_s/\bar{R}_b = 0.4$) without an electric field, the asymmetric response is marginal, resulting in negligible deviation from the principal path. 

Critically, under identical volumetric constraints, increased electrical loading reduces the internal pressure, underscoring the dominance of electrostatic softening over geometric stiffening in the post-buckling regime. 
Furthermore, the electric field modulates the spatial extent of the localised mode: higher voltages confine the deformation to a smaller region, intensifying the strain localisation.
\begin{figure}[t]
\centering
\begin{subfigure}[b]{0.8\linewidth} 
  \smash{\raisebox{0.7\height}{\makebox[0pt][l]{(a)}}}%
  \includegraphics[width=\linewidth]{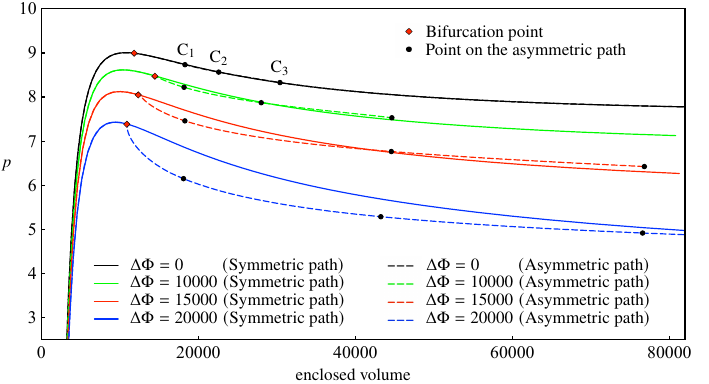}
  \label{fig:0.4_q}
\end{subfigure} 
\hfill
\begin{subfigure}[b]{0.8\linewidth}
  \smash{\raisebox{0.7\height}{\makebox[0pt][l]{(b)}}}%
  \includegraphics[width=\linewidth]{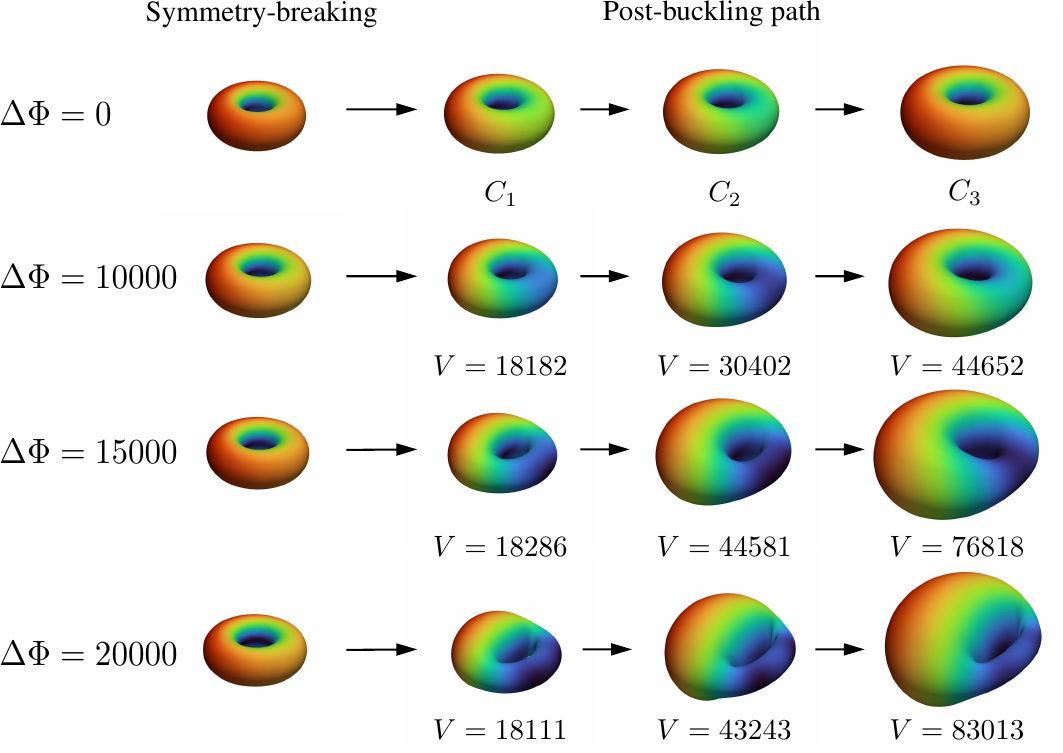}
  \label{fig:0.4_q_bx}
\end{subfigure}
\caption{Pressure–volume responses and deformed configurations illustrating the evolution of symmetry-breaking and post-buckling paths. Under increasing electric potential, the structure undergoes a transition from stable symmetric inflation to an asymmetric localized mode, with the maximum enclosed volume progressively increasing.  The torus has $\bar{R}_b = 10$ and $\bar{R}_s = 4$.}
\label{fig:0.4_q}
\end{figure}

\section{Conclusions}
\label{sec:conclusions}

This manuscript has presented a comprehensive isogeometric framework for the analysis of electroelastic thin shells, with a particular focus on capturing symmetry-breaking instabilities and post-buckling behaviour under strong electromechanical coupling. The key contributions are summarised as follows:

\begin{enumerate}
\item A nonlinear Kirchhoff--Love shell formulation for dielectric elastomers was developed, incorporating finite deformation kinematics, electromechanical coupling, and the plane stress condition. The total stress was decomposed into hyperelastic, Maxwell, and hydrostatic contributions, and the consistent tangent moduli were derived in closed form with static condensation to eliminate the thickness stress component.

\item The numerical implementation employs Catmull--Clark subdivision surfaces, providing the $C^1$-continuous basis functions required by Kirchhoff--Love shell theory, which enables accurate geometric representation and smooth deformation fields essential for capturing complex instability patterns.

\item By employing eigenmode perturbation to break symmetry, the evolution from axisymmetric to non-axisymmetric deformation modes was successfully captured, and the staged arc-length procedure allowed robust tracing of post-buckling equilibrium paths.

\item Numerical examples on spherical, prestretched circular plate, and toroidal membranes validated the framework's ability to capture large deformations, bifurcation onset, mode switching, and post-buckling responses under coupled electromechanical loading.
\end{enumerate}

Three numerical examples validated the proposed framework. The spherical membrane benchmark confirmed excellent agreement with the analytical solution across multiple potential differences, demonstrating the correct implementation of the electromechanical coupling. The prestretched circular plate example revealed the competitive interplay between mechanical prestress and electro-softening: higher prestretch suppresses the voltage-induced softening effect, while larger electric potentials significantly reduce the membrane's load-bearing capacity. The toroidal membrane example showcased the framework's unique capability to trace post-bifurcation equilibrium paths under electromechanical loading. The results demonstrated that electrical loading not only reduces the limit point pressure but also influences the post-buckling deformation pattern, with larger electric potentials producing more localised asymmetric deformations.

The proposed isogeometric framework provides a robust and accurate tool for analysing the complex nonlinear behaviour of dielectric elastomer thin shells. Its ability to handle large deformations, electromechanical coupling, and symmetry-breaking instabilities makes it well-suited for the design and optimisation of soft actuators, sensors, and energy harvesters. Future work will focus on extending the framework to incorporate dynamic effects, viscoelastic material behaviour, and multi-layer dielectric elastomer laminates.

\begin{nolinenumbers}
\section*{Acknowledgment}
Zhaowei Liu acknowledges the support from the National Natural Science Foundation of China (NSFC) under Grant No. 12502228.

\section*{Conflict of interest}
The authors declare that they have no conflict of interest.

\section*{Availability of data}
The datasets generated during the current study are available from the corresponding author upon reasonable request.

\section*{Availability of Code}
The code generated during the current study is available from the corresponding author on reasonable request.

\appendix
\section{Appendix: Derivation of the Invariant Form of Electric Energy}
\label{app:elec_energy_invariants}

For an incompressible electroelastic material, the electric energy density per unit reference volume is given by the linear dielectric model:
\begin{equation}
    \widetilde{W}_{\mathrm{elec}}(\mathbf{C},\bar{\mathbf{E}}) = -\frac{1}{2}\epsilon\,\bigl[\bar{\mathbf{E}}\otimes\bar{\mathbf{E}}\bigr]:\mathbf{C}^{-1},
    \label{eq:app_elec_definition}
\end{equation}
To express $W_{\mathrm{elec}}$ solely in terms of the scalar invariants, we eliminate $\mathbf{C}^{-1}$ using the Cayley--Hamilton theorem. For a $3\times3$ tensor $\mathbf{C}$, the theorem states:
\begin{equation}
    \mathbf{C}^3 - I_1\mathbf{C}^2 + I_2\mathbf{C} - I_3\mathbf{I} = \mathbf{0},
    \label{eq:app_cayley_hamilton}
\end{equation}
with $I_1 = \operatorname{tr}\mathbf{C}$, $I_2 = \frac{1}{2}\bigl[(\operatorname{tr}\mathbf{C})^2 - \operatorname{tr}(\mathbf{C}^2)\bigr]$, and $I_3 = \det\mathbf{C}$.
Under the incompressibility constraint, $I_3 = \det\mathbf{C} = J^2 = 1$. Substituting $I_3 = 1$ into Eq.~\eqref{eq:app_cayley_hamilton} and rearranging yields:
\begin{equation}
    \mathbf{C}^{-1} = \mathbf{C}^2 - I_1\mathbf{C} + I_2\mathbf{I}.
    \label{eq:app_C_inverse}
\end{equation}
Substituting Eq.~\eqref{eq:app_C_inverse} into Eq.~\eqref{eq:app_elec_definition} gives:
\begin{align}
    \bigl[\bar{\mathbf{E}}\otimes\bar{\mathbf{E}}\bigr]:\mathbf{C}^{-1}
    &= \bigl[\bar{\mathbf{E}}\otimes\bar{\mathbf{E}}\bigr]:\bigl[\mathbf{C}^2 - I_1\mathbf{C} + I_2\mathbf{I}\bigr] \notag \\
    &= \bar{\mathbf{E}}\cdot\mathbf{C}^2\bar{\mathbf{E}} - I_1\,\bar{\mathbf{E}}\cdot\mathbf{C}\bar{\mathbf{E}} + I_2\,\bar{\mathbf{E}}\cdot\bar{\mathbf{E}}.
    \label{eq:app_double_contraction}
\end{align}
Therefore, Eq.~\eqref{eq:app_double_contraction} becomes:
\begin{equation}
    \bigl(\bar{\mathbf{E}}\otimes\bar{\mathbf{E}}\bigr):\mathbf{C}^{-1} = I_6 - I_1 I_5 + I_2 I_4.
\end{equation}
Finally, substituting back into the electric energy expression yields the desired invariant form:
\begin{equation}
    \widetilde{W}_{\mathrm{elec}} = -\frac{1}{2}\epsilon\,\bigl[I_6 - I_1 I_5 + I_2 I_4\bigr].
    \label{eq:app_elec_final}
\end{equation}
\bibliographystyle{elsarticle-num-names}
\bibliography{sn-bibliography}
\end{nolinenumbers}
\end{document}